# COMPLEXITIES OF CONVEX COMBINATIONS AND BOUNDING THE GENERALIZATION ERROR IN CLASSIFICATION

By Vladimir Koltchinskii[1] and Dmitry Panchenko[2]

*University of New Mexico and Massachusetts Institute of Technology*

We introduce and study several measures of complexity of functions from the convex hull of a given base class. These complexity measures take into account the sparsity of the weights of a convex combination as well as certain clustering properties of the base functions involved in it. We prove new upper confidence bounds on the generalization error of ensemble (voting) classification algorithms that utilize the new complexity measures along with the empirical distributions of classification margins, providing a better explanation of generalization performance of large margin classification methods.

**1. Introduction.** Since the invention of *ensemble classification methods* (such as boosting), the convex hull conv($\mathcal{H}$) of a given base function class $\mathcal{H}$ has become an important object of study in the machine learning literature. The reason is that the ensemble algorithms typically output classifiers that are convex combinations of simple classifiers selected by the algorithm from the base class $\mathcal{H}$, and, because of this, measuring the complexity of the whole convex hull as well as of its subsets becomes very important in analysis of the generalization error of ensemble classifiers. Another important feature of boosting and many other ensemble methods is that they belong to the class of so-called *large margin methods*, that is, they are based on optimization of the empirical risk with respect to various loss functions that penalize not only for a misclassification (a negative classification margin), but also for a correct classification with too small positive margin. Thus, the very nature of these methods is to produce classifiers that tend to have rather large positive classification margins on the training data. Finding such classifiers becomes

Received May 2003; revised March 2004.
[1]Supported in part by NSF Grant DMS-03-04861 and NSA Grant MDA904-02-1-0075.
[2]Supported in part by an AT&T Buchsbaum Grant.

possible since the algorithms search for them in rather huge function classes (such as convex hulls of typical VC-classes used in classification).

This paper continues the line of research started by Schapire, Freund, Bartlett and Lee in [28] and further pursued in [2, 16, 19, 21, 26]. In these papers, the authors were trying to develop bounds on the generalization error of combined classifiers selected from the convex hull $\mathrm{conv}(\mathcal{H})$ in terms of the empirical distributions of their margins, as well as certain measures of complexity of the whole convex hull or its subsets to which the classifiers belong. Our main goal here is to suggest new margin type bounds that are based to a greater extent on complexity measures of *individual classifiers* from the convex hull. These bounds are more adaptive and more flexible than the previously known bounds (but they are also harder to prove). They take into account various properties of the convex combinations that are related to their generalization performance as classifiers, such as the sparsity of the weights and clustering properties of base functions.

The following notation and definitions will be used throughout the paper. Let $\mathcal{X}$ be a measurable space (space of instances) and let $\mathcal{Y} = \{-1, +1\}$ be the set of labels. Let $\mathbb{P}$ be a probability measure on $\mathcal{X} \times \mathcal{Y}$ that describes the underlying distribution of instances and their labels. We do not assume that the label $y$ is a deterministic function of $x$; in general, it can also be random, which means that the conditional probability $\mathbb{P}(y = 1|x)$ may be different from 0 or 1. Let $\mathcal{H}$ be a class of measurable functions $h : \mathcal{X} \to [-1, 1]$. Denote by $\mathcal{P}(\mathcal{H})$ the set of all discrete distributions on $\mathcal{H}$ and let $\mathcal{F}$ be the convex hull of $\mathcal{H}$,

$$\mathcal{F} = \mathrm{conv}(\mathcal{H}) := \left\{ \int h(\cdot) \lambda(dh) : \lambda \in \mathcal{P}(\mathcal{H}) \right\}.$$

For $f \in \mathcal{F}$ we assume that $\mathrm{sign}(f(x))$ is used to classify $x \in \mathcal{X}$ [$\mathrm{sign}(f(x)) = 0$ meaning that no decision is made]. Functions $f \in \mathcal{F}$ are sometimes called *voting classifiers*, since for a convex combination $f = \sum \lambda_j h_j$ the weight (coefficient) $\lambda_j$ can be interpreted as the voting power of an individual classifier $h_j$ (they are also called *ensemble classifiers*). The generalization error of any classifier $f \in \mathcal{F}$ is defined as

(1.1) $$\mathbb{P}(\mathrm{sign}(f(x)) \neq y) = \mathbb{P}(yf(x) \leq 0).$$

Given an i.i.d. sample $(X_1, Y_1), \ldots, (X_n, Y_n)$ from the distribution $\mathbb{P}$, let $\mathbb{P}_n$ denote its *empirical distribution*. For a measurable function $g$ on $\mathcal{X} \times \mathcal{Y}$, denote

$$\mathbb{P}g = \int g(x, y) \, d\mathbb{P}(x, y), \qquad \mathbb{P}_n g = n^{-1} \sum_{i=1}^{n} g(X_i, Y_i).$$

Whenever it is needed, we use the same notation $\mathbb{P}g, \mathbb{P}_n g$ or $\mathbb{P}(A), \mathbb{P}_n(A)$ for functions $g$ that depend only on $x$ and for sets $A \subset \mathcal{X}$ (the meaning of the notation in this case is obvious). The probability measure on the main sample space (on which all the random variables including the training examples are defined) will be denoted by **Pr** (not to confuse it with $\mathbb{P}$).

In the paper we study the generalization error (1.1) of classifiers from the convex hull of a class $\mathcal{H}$ which is typically assumed to be "small," a condition that is described precisely in terms of some complexity assumptions on $\mathcal{H}$ [see (2.2)]. A number of popular classification algorithms output classifiers of this type. Below we briefly discuss two of them: *AdaBoost*, which is the most well-known classification algorithm of *boosting* type, and also *bagging*. We provide some heuristic explanations of why these algorithms might have a tendency to output convex combinations of classifiers from the base class with a certain degree of sparsity of their weights and clustering of the base classifiers.

*AdaBoost.* The algorithm starts by assigning equal weights $w_j^{(1)} = \frac{1}{n}$ to all the training examples $(X_j, Y_j)$. At iteration number $k$, $k = 1, \ldots, T$, the algorithm attempts to minimize the weighted training error with weights $w_j^{(k)}$ over the base class $\mathcal{H}$ of functions $h : S \mapsto \{-1, 1\}$ (such that $h \in \mathcal{H}$ implies $-h \in \mathcal{H}$). If $e_k$ denotes the weighted training error of the approximate solution $h_k$ of this minimization problem, the algorithm computes the coefficient

$$\alpha_k := \frac{1}{2} \log \frac{1 - e_k}{e_k},$$

which is nonnegative since $e_k \leq \frac{1}{2}$, and then updates the weights according to the formula

$$w_j^{(k+1)} := \frac{w_j^{(k)} e^{-Y_j \alpha_k h_k(X_j)}}{Z},$$

where $Z$ is a normalizing constant that makes the weights add up to 1. After $T$ iterations, the algorithm outputs the classifier $f = \sum_{k=1}^{T} \lambda_k h_k$, where

$$\lambda_k := \frac{\alpha_k}{\sum_{j=1}^{T} \alpha_j}.$$

Typically, the class $\mathcal{H}$ is relatively small so that it is easy to design an efficient algorithm (often called a weak learner) of approximate minimization of the weighted training error over the class. The result of this, however, is that at many iterations the weak learner outputs classifiers $h_k$ from the base $\mathcal{H}$ whose weighted training error is just a little smaller than $1/2$. If this is the case at iteration $k$, the coefficient $\alpha_k$ is close to 0 and the weights $w_j^{(k+1)}$ do not differ much from the weights $w_j^{(k)}$. If the weak learner possesses some stability, this means that the base classifier $h_{k+1}$ is close to the base classifier $h_k$. As a result, when the algorithm proceeds one observes a slow drift of the classifiers $h_k$ in the "hypotheses space" $\mathcal{H}$, and the coefficients of these classifiers in the convex combination will be small until we reach a place in $\mathcal{H}$ where the stability of the weak learner breaks down and it outputs a classifier with a weighted training error significantly smaller than $1/2$. Thus, one can expect a certain degree of sparsity (many small coefficients) and of clustering (many base classifiers that are close to one another) of the resulting convex combination.

*Bagging* [9]. The algorithm at each iteration produces a bootstrap sample drawn from the training data and outputs a classifier that minimizes the corresponding bootstrap training error over the base class $\mathcal{H}$. After $T$ iterations the algorithm averages the resulting $T$ base classifiers, creating a convex combination with equal weights $\lambda_k := \frac{1}{T}$. Again, if the weak learner possesses some stability and since each bootstrap sample is a "small perturbation" of the training data, one can expect some degree of clustering of the base functions involved in the convex combination. (In this case, it is impossible to talk about the sparsity of the coefficients since all of them are equal.)

These explanations are of course rather heuristic in nature and somewhat vague. The reality might be much more complicated since, for instance, weak learners are not necessarily stable. Often, lack of stability of the weak learner is viewed as an advantage since it allows the algorithm to create more "diverse" ensembles of base classifiers and to produce a combined classifier with larger margins. However, the bounds of this paper seem to suggest that the performance of combined classifiers is related to a rather delicate trade-off between their complexity and margin properties. So, stability of the weak learner is a good and a bad property at the same time (one should rather talk about optimal stability). The phenomenon of sparsity of the coefficients is much better understood in the case of support vector machines (see [30] for recent results in this direction) and the development of these ideas for ensemble methods remains an open problem that is beyond the scope of our paper. However, regardless of how close this explanation is to the truth, some degree of sparsity and clustering in convex combinations output by popular learning algorithms can be observed in experiments (see some very preliminary results in [20] and more results in [1]). Our intention here is not to study why this is happening, but rather to understand what kind of influence sparsity and clustering properties of convex combinations output by *AdaBoost* and other classification algorithms have on their generalization performance.

Another motivation to study the complexities based on sparsity and clustering comes from learning theory, where it has become common to use global or localized complexities based on sup-norm or continuity modulus of empirical or Rademacher processes involved in the problem and indexed by the class $\mathcal{F}$ in order to bound the generalization error (see [5, 8, 17, 18, 23]). However, these complexities do not necessarily measure the accuracy of modern classification methods correctly. The reason is that they are based on deviations of the empirical measure $\mathbb{P}_n$ from the true distribution $\mathbb{P}$ uniformly over the whole class $\mathcal{F}$ or over $\mathcal{L}_2(\mathbb{P})$-balls in the class, while the learning algorithms might have some intrinsic ways to restrict complexities of the classifiers they output by searching for a minimum of empirical risk in some parts of the class $\mathcal{F}$ with restricted complexity (although this part is typically data-dependent, cannot be specified in advance and has to be determined in a rather complicated model selection process). Thus, there is a need to develop new more adaptive bounds that take into account complexities of individual classifiers in the class and can be applied to the classifiers output by learning algorithms. A possible general approach to such complexities can be described as follows. Suppose $\{\mathcal{G}\}$ is a family of subclasses of the class $\mathcal{F}$ and let $c_n(\mathcal{G})$ be a complexity measure associated with the class $\mathcal{G}$ (e.g., it can be based on a localized Rademacher complexity of $\mathcal{G}$). Suppose also it has been observed that a learning algorithm tends to output classifiers from subclasses $\mathcal{G}$ with small values of complexity $c_n(\mathcal{G})$ ("sparse subclasses"). Then a natural question to ask is whether the quantity of the type $c_n(f) := \inf\{c_n(\mathcal{G}) : \mathcal{G} \ni f\}$ (which is already an individual complexity of $f$) provides bounds on the generalization error of $f$. In the case where $\{\mathcal{G}\}$ is a countable family of nested subclasses, such questions are related to structural risk minimization and other model selection techniques. However, in classification one often encounters more complicated situations, such as the setting of Theorem 5 below, where a natural family $\{\mathcal{G}\}$ is neither countable nor nested and consists of distribution-dependent classes indexed by a functional parameter (see the definition of the classes $\mathcal{F}_{Q,p,N}^C$ before Lemma 2). The study of complexity measures that occur in such more complicated model selection frameworks is our main subject here. In the next section we will try to develop several new approaches to measuring complexities of convex combinations and use these complexities in new bounds on generalization error in classification.

**2. Main results.** The first important result about the generalization error of classifiers from $\mathcal{F} = \text{conv}(\mathcal{H})$ was proved in [28], where the generalization ability of voting classifiers is explained in terms of the empirical distribution $\mathbb{P}_n(yf(x) \leq \delta)$ of the quantity $yf(x)$ called *margin*. The authors prove that if $\mathcal{H} = \{2I(x \in C) - 1 : C \in \mathcal{C}\}$, where $\mathcal{C}$ is a Vapnik–Chervonenkis class of sets with VC-dimension $V$ (for definitions see, e.g., [32] or [12]), then for all $t > 0$ with probability at least $1 - e^{-t}$ for all $f \in \mathcal{F} = \text{conv}(\mathcal{H})$ we have

$$\mathbb{P}(yf(x) \leq 0)$$

(2.1)
$$\leq \inf_{\delta \in (0,1]} \left( \mathbb{P}_n(yf(x) \leq \delta) + K\left( \left( \frac{V \log^2(n/\delta)}{n\delta^2} \right)^{1/2} + \left( \frac{t}{n} \right)^{1/2} \right) \right),$$

where $K > 0$ is an absolute constant. To understand this result, let us give one interpretation of the margin $yf(x)$. One can think of $yf(x)$ as the "confidence" of prediction of the example $x$, since $f$ classifies $x$ correctly if and only if $yf(x) > 0$; and if $f(x)$ is large in absolute value it means that it makes its prediction with high confidence. If $f$ classifies most of the training examples with high confidence, then for some $\delta > 0$ (which is not "too small") the proportion of examples $\mathbb{P}_n(yf(x) \leq \delta)$ classified below the confidence $\delta$ will be small. The second term of the bound is of the order $(\sqrt{n}\delta)^{-1}$, and will also be small for large $n$, which makes the bound meaningful.

This result was extended by Schapire and Singer in [29] to classes of real-valued functions, namely, to so-called VC-subgraph classes (for definition see [32]), and was further extended in several directions in [19] and [21]. The main idea of this follow-up work was to replace the second term of the bound proved by Schapire et al. [28] by a function $\varepsilon_n(\mathcal{F}; \delta; t)$ that has better dependence on the sample size $n$ and on the margin parameter $\delta$. The bounds obtained in [19] are also more general: they apply to arbitrary function classes $\mathcal{F}$, not only to the convex hulls.

Given a probability distribution $Q$ on $\mathcal{X}$ and a class $\mathcal{H}$ of measurable functions on $\mathcal{X}$, denote
$$d_{Q,2}(f,g) := (Q(f-g)^2)^{1/2}, \qquad f, g \in \mathcal{H},$$
the $\mathcal{L}_2(Q)$-distance in $\mathcal{H}$. Let the *covering number* $N_{d_{Q,2}}(\mathcal{H}, u)$ be the minimal number of $d_{Q,2}$-balls of radius $u > 0$ with centers in $\mathcal{H}$ needed to cover $\mathcal{H}$. The logarithm of this number $H_{d_{Q,2}}(\mathcal{H}, u) := \log N_{d_{Q,2}}(\mathcal{H}, u)$ is called the $u$-entropy of $\mathcal{H}$ with respect to $d_{Q,2}$. In what follows, we will also use $\mathcal{L}_p(Q)$-distances and the corresponding covering numbers and entropies for $p \in [1, +\infty]$.

Often, it makes sense to assume (and it will be assumed in what follows) that the family of weak classifiers $\mathcal{H}$ satisfies the condition

(2.2)
$$\sup_{Q \in \mathcal{P}(\mathcal{X})} N_{d_{Q,2}}(\mathcal{H}, u) = \mathcal{O}(u^{-V})$$

for some $V > 0$, where $\mathcal{P}(\mathcal{X})$ is the set of all discrete distributions on $\mathcal{X}$. For example, if $\mathcal{H}$ is a VC-subgraph class with VC-dimension $V(\mathcal{H})$, then by the well-known result that goes back to Dudley and Pollard (see [14] for the current version), (2.2) holds with $V = 2V(\mathcal{H})$, namely,

(2.3)
$$\sup_{Q \in \mathcal{P}(\mathcal{X})} N_{d_{Q,2}}(\mathcal{H}, u) \leq e(V(\mathcal{H})+1) \left( \frac{2e}{u^2} \right)^{V(\mathcal{H})}.$$

Under the condition (2.2), the bound (2.1) was slightly improved by Koltchinskii and Panchenko in [19], who proved that for all $t > 0$ with probability at least $1 - e^{-t}$ for all $f \in \mathcal{F} = \text{conv}(\mathcal{H})$ we have

$$\mathbb{P}(yf(x) \leq 0) \leq \inf_{\delta \in (0,1]} \left( \mathbb{P}_n(yf(x) \leq \delta) + K\left( \left(\frac{V}{n\delta^2}\right)^{1/2} + \left(\frac{t}{n}\right)^{1/2} \right) \right),$$

thus getting rid of the logarithmic factor $\log^2(n/\delta)$ in the second term of (2.1). By itself this improvement is insignificant, but the generality of the methods developed in [19] allowed the authors to obtain this type of bound for general classes $\mathcal{F}$ of classifiers (not necessarily the convex hulls of VC-classes) and to make some significant improvements in other situations, for example, for neural networks. (The first margin type bounds for general function classes, including neural networks, were based on $\mathcal{L}_\infty$-entropies and shattering dimensions of the class; see [4].) Moreover, it was shown in [19] that (2.1) can be further improved in the so-called zero-error case, when $\mathbb{P}_n(yf(x) \leq \delta)$ is small for $\delta \to 0$. Namely, the following result holds. Assume that $\mathcal{H}$ satisfies (2.2) and let $\alpha = 2V/(V+2)$. Then, for all $t > 0$ with probability at least $1 - e^{-t}$ for all $f \in \mathcal{F}$ we have (with some numerical constant $K > 0$)

$$(2.4) \quad \begin{aligned} &\mathbb{P}(yf(x) \leq 0) \\ &\leq K \inf_{\delta \in (0,1]} \left( \mathbb{P}_n(yf(x) \leq \delta) + \left( \left(\frac{1}{\delta}\right)^{2\alpha/(2+\alpha)} n^{-2/(\alpha+2)} + \frac{t}{n} \right) \right). \end{aligned}$$

This bound will be meaningful if

$$\delta^* = \sup\{\delta : \delta^{2\alpha/(2+\alpha)} \mathbb{P}_n(yf(x) \leq \delta) \leq n^{-2/(2+\alpha)}\}$$

is not "too small," which means that $\mathbb{P}_n(yf(x) \leq \delta)$ should decrease "fast enough" when $\delta \to 0$. Actually, this bound holds not only for classes of functions $\mathcal{F} = \text{conv}(\mathcal{H})$ where $\mathcal{H}$ satisfies (2.2), but for any class $\mathcal{F}$ such that

$$(2.5) \quad \sup_{Q \in \mathcal{P}(\mathcal{X})} \log N_{d_{Q,2}}(\mathcal{F}, u) = \mathcal{O}(u^{-\alpha}), \qquad \alpha \in (0, 2),$$

or even when the uniform entropy in (2.5) is replaced by the entropy with respect to empirical $\mathcal{L}_2$-distance $d_{\mathbb{P}_n,2}$. It is well known that the convex hull $\mathcal{F} = \text{conv}(\mathcal{H})$ of the class $\mathcal{H}$ satisfying (2.2) satisfies (2.5) with $\alpha = 2V/(V+2)$ (see, e.g., [32]), which explains a particular choice of $\alpha$ in (2.4). Under the condition (2.5) on $\mathcal{F}$ the bound of (2.4) is optimal as shown in [19] by constructing a special class of functions $\mathcal{F}$ in Banach space $l_\infty$ of uniformly bounded sequences. Finally, note that the constant $K$ involved in the bound can be redistributed between the two terms: in front of the term $\mathbb{P}_n(yf(x) \leq \delta)$ one can put a constant arbitrarily close to 1 at the price of making the constant in front of the second term large.

In [21] Koltchinskii, Panchenko and Lozano proved the bounds on generalization error under more general assumption on the entropy of the class $\mathcal{F}$:

$$(2.6) \qquad \int_0^x H_{d_{\mathbb{P}_n,2}}^{1/2}(\mathcal{F}; u)\, du \leq D\psi(x), \qquad x > 0,$$

with some constant $D > 0$ and with a concave function $\psi$. They showed that in this case the term

$$\left(\frac{1}{\delta}\right)^{2\alpha/(2+\alpha)} n^{-2/(\alpha+2)}$$

involved in the bound (2.4) should be replaced by the quantity $\varepsilon_n^\psi(\delta)$ defined as the largest solution of the equation

$$\varepsilon = \frac{1}{\delta\sqrt{n}}\psi(\delta\sqrt{\varepsilon}),$$

leading to so-called $\psi$-bounds on generalization error.

Margin-type bounds on generalization error can be also expressed in terms of other entropies, in particular, $\mathcal{L}_\infty$-entropy and in terms of shattering dimension of the class, as in the papers of Bartlett [4] (that preceded [28]) and of Antos, Kégl, Linder and Lugosi [2]. A typical bound in terms of $\mathcal{L}_\infty$-entropy is of the form

$$(2.7) \quad \mathbb{P}(yf(x) \leq 0) \leq K \inf_{\delta \in (0,1]} \left(\mathbb{P}_n(yf(x) \leq \delta) + \frac{\log \mathbb{E} N_{d_{\mathbb{P}_n,\infty}}(\mathcal{F}; \delta/2) + t}{n}\right)$$

for all $f \in \mathcal{F}$ with probability at least $1 - e^{-t}$. The $\mathcal{L}_\infty$-entropy is always larger than $\mathcal{L}_2$-entropy, but for special classes of functions the difference might be not very significant, and because of a different form the $\mathcal{L}_\infty$-bound has sometimes an advantage over the $\mathcal{L}_2$-bounds. However, the detailed comparison of these bounds goes beyond the scope of this paper.

Numerous experiments with *AdaBoost* and some other classification algorithms showed that in practice the bounds of type (2.4) hold with smaller values of $\alpha$ than the theoretical considerations (based on the estimates of the entropy of the *whole* convex hull) suggest. This means that ensemble classifiers often belong to a subset of the convex hull of a smaller entropy than the entropy of the whole convex hull. A natural question is whether it is possible to incorporate in the bound on generalization error the information about the *individual* complexity of the actual classifier rather than use *global* complexity of the whole convex hull. In other words, is it possible to replace the function $\psi$ from condition (2.6) by a data-dependent and classifier-dependent function that would make the $\psi$-bounds on generalization error more adaptive?

The fact that the margin type bounds hold in such generality means, at least on the intuitive level, that the explicit structure of the convex hull is not used there. On the contrary, in this paper we will heavily utilize the structure of the convex hull and prove new bounds that reflect some measures of complexity of convex combinations.

The idea of using a certain measure of complexity of *individual* convex combinations already appeared in [21], where the authors suggested a way to use a rate of decay of weights $\lambda_j$ in the convex combination $f = \sum_{j=1}^T \lambda_j h_j$ to improve the bound on the generalization error of $f$. This measure, called approximate $\gamma$-dimension, is defined as follows. Let us assume that the weights are arranged in the decreasing order $|\lambda_1| \geq |\lambda_2| \geq \cdots$. For a number $\gamma \in [0,1]$, *the approximate $\gamma$-dimension* of $f$ is defined as the smallest integer number $d \geq 0$ such that there exist $T \geq 1$, functions $h_j \in \mathcal{H}$, $j = 1, \ldots, T$, and numbers $\lambda_j \in \mathbb{R}$, $j = 1, \ldots, T$, satisfying the conditions $f = \sum_{j=1}^T \lambda_j h_j$, $\sum_{j=1}^T |\lambda_j| \leq 1$ and $\sum_{j=d+1}^T |\lambda_j| \leq \gamma$. Note that in [21] the authors dealt with the symmetric convex hull, so the coefficients $\lambda_j$ are not necessarily positive. The $\gamma$-dimension of $f$ will be denoted by $d(f; \gamma)$.

Then, for all $t > 0$ with probability at least $1 - e^{-t}$ we have for all $f \in \mathcal{F} = \mathrm{conv}(\mathcal{H})$ (again with $\alpha = \frac{2V}{V+2}$)

$$\mathbb{P}(yf(x) \leq 0)$$

$$(2.8) \quad \leq K \inf_{\delta \in (0,1]} \left( \mathbb{P}_n(yf(x) \leq \delta) \right.$$

$$\left. + \inf_\gamma \left( \frac{d(f;\gamma)}{n} \log \frac{n}{\delta} + \left(\frac{\gamma}{\delta}\right)^{2\alpha/(2+\alpha)} n^{-2/(\alpha+2)} + \frac{t}{n} \right) \right).$$

This is an improvement over (2.4), which can be seen by comparing the infimum over $\gamma$ of the expression in the bound with the value of the expression for $\gamma = 1$ and noting that $d(f; 1) = 0$. For example, if the weights decrease polynomially $|\lambda_j| \sim j^{-\beta}, \beta > 1$, or exponentially $|\lambda_j| \sim e^{-\beta j}, \beta > 0$, then explicit minimization over $\gamma$ shows that in these cases (2.8) can be a substantial improvement over (2.4) (see examples in [21]).

Our first result in this paper also deals with bounding the generalization error of a classifier $f = \sum_{j=1}^T \lambda_j h_j \in \mathcal{F} = \mathrm{conv}(\mathcal{H})$ in terms of complexity measures taking into account the sparsity of the weights $\lambda_j$. Theorem 1 below is a new version of the results of [21] [specifically, of the bound (2.8)] that can be interpreted as interpolation between zero-error and nonzero-error cases; as its corollary we will give a new short proof of (2.8). Theorem 2 is another result in this direction with a different dependence of the bound on the sample size and the margin parameter $\delta$.

Let $\Phi = \{\varphi_\delta : \mathbb{R} \to [0,1] : \delta \in \Delta \subset \mathbb{R}_+\}$ be a countable family of Lipschitz functions such that the Lipschitz norm of $\varphi_\delta$ is bounded by $\delta^{-1}$, that is,

$$|\varphi_\delta(s_1) - \varphi_\delta(s_2)| \leq \delta^{-1}|s_1 - s_2|,$$

and $\sum_{\delta \in \Delta} \delta < \infty$. In applications, such functions are frequently used as loss functions in empirical risk minimization procedures of boosting type that output large margin classifiers. One can use a specific choice of $\Delta = \{2^{-k} : k \geq 1\}$. The following theorem holds.

THEOREM 1. *If* (2.2) *holds, then for all* $t > 0$ *with probability at least* $1 - e^{-t}$ *for all* $f \in \mathcal{F} = \mathrm{conv}(\mathcal{H})$ *and* $\delta \in \Delta = \{2^{-k} : k \geq 1\}$,

$$\frac{\mathbb{P}\varphi_\delta(yf(x)) - \mathbb{P}_n\varphi_\delta(yf(x))}{(\mathbb{P}\varphi_\delta(yf(x)))^{1/2}}$$

$$\leq K \inf_\gamma \left( \left(\frac{d(f;\gamma)}{n} \log \frac{n}{\delta}\right)^{1/2} + \left(\frac{\gamma}{\delta}\right)^{\alpha/2} \frac{(\mathbb{P}\varphi_\delta(yf(x)))^{-\alpha/4}}{n^{1/2}} + \left(\frac{t}{n}\right)^{1/2} \right),$$

*where* $\alpha = 2V/(V+2)$.

Let us take, for example, $\varphi_\delta$ such that $\varphi_\delta(s) = 1$ for $s \leq 0$, $\varphi_\delta = 0$ for $s \geq \delta$ and $\varphi_\delta$ is linear for $0 \leq s \leq \delta$. For any probability measure $Q$ (e.g., $Q = \mathbb{P}$ or $\mathbb{P}_n$), one can write

(2.9) $$Q(yf(x) \leq 0) \leq Q\varphi_\delta(yf(x)) \leq Q(yf(x) \leq \delta).$$

For this choice of $\varphi_\delta$ and for a fixed $f$ let us denote $a = \mathbb{P}\varphi_\delta(yf(x))$ and $b = \mathbb{P}_n\varphi_\delta(yf(x))$. It is clear that after minimizing the expression involved in the right-hand side over $\gamma$, the inequality of Theorem 1 can be written as

$$a \leq b + ua^{1/2} + va^{1/2-\alpha/4},$$

where $u$ and $v$ are constants depending on the parameters involved in the inequality. Since the right-hand side of the last inequality is strictly concave with respect to $a$, this inequality can be uniquely solved for $a$ or, in other words, it can be equivalently written as $a \leq \rho(b)$ for unique positive function $\rho$, which is, obviously, increasing in $b$. Combining this with (2.9) we get

$$\mathbb{P}(yf(x) \leq 0) \leq \mathbb{P}\varphi_\delta(yf(x)) \leq \rho(\mathbb{P}_n\varphi_\delta(yf(x))) \leq \rho(\mathbb{P}_n(yf(x) \leq \delta)).$$

The analysis of $\rho$ will readily imply the main result in [21].

COROLLARY 1. *If* (2.2) *holds and* $\alpha = 2V/(V+2)$, *then for any* $t > 0$ *with probability at least* $1 - e^{-t}$ (2.8) *holds for all* $f \in \mathcal{F} = \mathrm{conv}(\mathcal{H})$.

Roughly speaking, Corollary 1 describes the zero-error case of Theorem 1. Thus, Theorem 1 is a more general and flexible formulation of the main result in [21], as it interpolates between zero- and nonzero-error cases.

Next we will present a new bound on the generalization error of voting classifiers that takes into account the sparsity of weights in the convex combination. Given $\lambda \in \mathcal{P}(\mathcal{H})$ and $f(x) = \int h(x)\lambda(dh)$, we can also represent $f$ as $f = \sum_{k=1}^{T} \lambda_k h_k$ with $T \leq \infty$ (since $\lambda$ is a discrete probability measure). Without loss of generality let us assume that $\lambda_1 \geq \lambda_2 \geq \cdots$. We define $\gamma_d(f) = \sum_{k=d+1}^{T} \lambda_k$ and for $\delta > 0$ we define the *effective dimension* function by

$$(2.10) \qquad e_n(f, \delta) = \min_{0 \leq d \leq T} \left( d + \frac{2\gamma_d^2(f)}{\delta^2} \log n \right).$$

This name is motivated by the fact that (as will become clear from the proof of Theorem 2 below) it can be interpreted as a dimension of a subset of the convex hull $\text{conv}(\mathcal{H})$ that contains a "good" approximation of $f$.

THEOREM 2 (Sparsity bound). *If* (2.2) *holds, then there exists an absolute constant $K > 0$ such that for all $t > 0$ with probability at least $1 - e^{-t}$ for all $\lambda \in \mathcal{P}(\mathcal{H})$ and $f(x) = \int h(x)\lambda(dh)$,*

$$\mathbb{P}(yf(x) \leq 0) \leq \inf_{\delta \in (0,1]} (U^{1/2} + (\mathbb{P}_n(yf(x) \leq \delta) + U)^{1/2})^2,$$

*where*

$$U = K\left( \frac{V e_n(f, \delta)}{n} \log \frac{n}{\delta} + \frac{t}{n} \right).$$

It follows from the bound of the theorem that for all $\varepsilon > 0$

$$\mathbb{P}(yf(x) \leq 0) \leq \inf_{\delta \in (0,1]} \left( (1+\varepsilon)\mathbb{P}_n(yf(x) \leq \delta) + \left(2 + \frac{1}{\varepsilon}\right) U \right),$$

which is a more explicit version of the result. Results of similar flavor can be, in principle, also obtained as a consequence of entropy-based margin-type bounds, in particular, using the $\mathcal{L}_\infty$-entropy. However, we believe that the more direct probabilistic argument we use in our proof (that goes back to [28]) is very natural in this problem. Moreover, the same argument is typically present in the derivation of entropy bounds for the convex hull or its subsets needed in alternative proofs. Taking this into account, the direct proof we give here is shorter and easier. This becomes especially clear in Theorems 3 and 4, where the entropy bounds on subsets of the convex hull with restrictions on the variance of convex combinations (see the definitions below) are most likely unknown. It is also worth mentioning that the same randomization idea combined with a couple of other techniques can be used in some other situations where probabilistic interpretation is not straightforward, for instance, for kernel machines and their hierarchies (see [1]).

The following result was proved in [11]. Let $\mathcal{H}$ be a finite class with $N = \text{card}(\mathcal{H})$ and let $\delta_*$ be the minimal margin on the training examples, that is,

$$\delta_* = \delta_*(f) = \min_{i \leq n} Y_i f(X_i) = \sup\{\delta : \mathbb{P}_n(yf(x) \leq \delta) = 0\}.$$

Then for any $t > 0$ with probability at least $1 - e^{-t}$ we have, for all $f \in \mathcal{F} = \text{conv}(\mathcal{H})$ such that $\delta_*(f) \geq (32/N)^{1/2}$,

(2.11)  $$\mathbb{P}(yf(x) \leq 0) \leq K\left(\frac{\log N}{n\delta_*^2} + \frac{t}{n}\right).$$

We notice that

$$e_n(f, \delta) = \min_{0 \leq d \leq T}\left(d + \frac{2\gamma_d^2(f)}{\delta^2}\log n\right) \leq \frac{2}{\delta^2}\log n,$$

where the last inequality follows by taking $d = 0$ in the expression under the infimum. This shows that as a corollary of Theorem 2 one can extend the result of Breiman [11] to much more general classes of functions [the role of $\log N$ in (2.11) being now played by $V \log n$]. Moreover, the bound of Theorem 2 interpolates between zero-error and nonzero-error cases without any assumptions on the empirical distribution of the margin $\mathbb{P}_n(yf(x) \leq \delta)$. To illustrate the role of the effective dimension $e_n(f, \delta)$ let us suppose that the weights $\lambda_j$ decrease polynomially or exponentially fast:

EXAMPLE.  (a) If $\lambda_j \sim j^{-\beta}$ for $\beta > 1$, then one can explicitly minimize the expression in (2.10), which in the zero-error case $\mathbb{P}_n(yf(x) \leq \delta_*) = 0$ gives

$$\mathbb{P}(yf(x) \leq 0) \leq K(\beta)\left(\frac{V}{n\delta_*^{2/(2\beta-1)}}\log^2\frac{n}{\delta_*} + \frac{t}{n}\right),$$

which can be a significant improvement for large values of $\beta$.

(b) If $\lambda_j \sim e^{-j}$, then again one can explicitly minimize the expression in (2.10), which in the zero-error case $\mathbb{P}_n(yf(x) \leq \delta_*) = 0$ gives

$$\mathbb{P}(yf(x) \leq 0) \leq K\left(\frac{V}{n}\log^2\frac{n}{\delta_*} + \frac{t}{n}\right).$$

It is quite clear that one can come up with many alternative definitions of sparsity measures of convex combinations that are based only on the sizes of coefficients. For instance, one can measure the size of the "tail" of the convex combination (after the $d$ largest coefficients have been removed) using a different norm instead of the $\ell_1$-norm we used above. However, our approach seems to be reasonable since it is based on the idea of splitting the whole convex combination into two parts, one of them being $d$-dimensional and another one belonging to a rescaled convex hull of $\mathcal{H}$ (the whole convex hull times a small coefficient, which is a natural "neighborhood" of 0 in the convex hull).

The major drawback of this type of bound, however, is that it takes into account only the size of the coefficients of the convex combination, but not the "closeness" of the base functions involved in it. Such a "closeness" (reflected, e.g., in the fact that the base functions classify most of the examples the same way or, more generally, can be divided into several groups with the functions within each group classifying similarly) could possibly lead to further complexity reduction.

We suggest below two different approaches to this problem. The first approach is based on interpreting the convex combination as a mean of a function $h$ randomly drawn from the class $\mathcal{H}$ with some probability distribution $\lambda$. Then in order to measure the complexity of the convex combination it becomes natural to bring in probabilistic quantities such as the variance of the convex combination introduced below. In the extreme case, when all classifiers $h_j$ are equal, $f$ belongs to a simple class $\mathcal{H}$ itself rather than to the possibly very large class $\mathcal{F}$; in this case, the variance is equal to 0 and this is reflected in our generalization analysis of $f$. This approach is clearly related to the randomization proof of margin type bounds in [28], but its real roots are in the well-known work of B. Maurey (see [27]) that provided a probabilistic argument often used in bounding the entropy of the convex hull. The approach might be also of interest to practitioners since variance can be easily incorporated in risk minimization techniques as a complexity penalty. The generalization bounds based on the notion of variance are given in Theorems 3 and 4.

The second approach does not rely on the probabilistic interpretation, but rather exploits the nonuniqueness of representing functions by convex combinations and is based on covering numbers of the set of base functions in "optimal" representations of $f$. Thus, the metric structure of the base class replaces in this approach the probabilistic structure. The generalization bound based on this approach is given in Theorem 5.

Despite the fact that, possibly, there might be many other ways to define complexities of this type, we believe that the approaches we are using have very natural connections to important mathematical structures involved in the problem.

Given $\lambda \in \mathcal{P}(\mathcal{H})$, consider

$$f(x) = \int h(x)\lambda(dh) = \sum_{k=1}^{T} \lambda_k h_k(x).$$

We ask the following question: what if the functions $h_1, \ldots, h_T$ are, in some sense, close to each other? For example, $n^{-1} \sum_{k=1}^n (h_i(X_k) - h_j(X_k))^2$ is small for all pairs $i, j$. In this case, the convex combination can be approximated "well" by only one function from $\mathcal{H}$. Or, more generally, one can imagine the situation when there are several clusters of functions among $h_1, \ldots, h_T$ such that within each cluster all functions are close to each other. This information should be reflected in the generalization error of classifier $f$, since it can be approximated by a classifier from a small subset of $\mathcal{F}$. Below we prove two results in this direction. We will start by describing the result where we consider $h_1, \ldots, h_T$ as one (hopefully "small") cluster, and then we will naturally generalize it to any number of clusters.

We define a pointwise variance of $h$ with respect to the distribution $\lambda$ by

$$(2.12) \qquad \sigma_\lambda^2(x) = \int \left( h(x) - \int h(x) \lambda(dh) \right)^2 \lambda(dh).$$

Clearly, $\sigma_\lambda^2(x) = 0$ if and only if

$$h(x) = \int h(x) \lambda(dh), \qquad \lambda\text{-a.e. on } \mathcal{H},$$

or, equivalently (in the case of a discrete measure $\lambda$), if $h_1(x) = h_2(x)$ for all $h_1, h_2 \in \mathcal{H}$ with $\lambda(\{h_1\}) > 0$, $\lambda(\{h_2\}) > 0$. The complexity characteristics of a similar flavor are sometimes used in the current work on PAC Bayesian bounds on generalization performance of aggregated estimates for least square regression; see [3].

THEOREM 3. *If* (2.2) *holds, then there exists an absolute constant $K > 0$ such that for all $t > 0$ with probability at least $1 - e^{-t}$ for all $\lambda \in \mathcal{P}(\mathcal{H})$ and $f(x) = f_\lambda(x) = \int h(x) \lambda(dh)$,*

$$\mathbb{P}(y f_\lambda(x) \leq 0)$$
$$\leq K \inf_{0 < \delta \leq \gamma \leq 1} \left( \mathbb{P}_n(y f_\lambda(x) \leq \delta) + \mathbb{P}_n(\sigma_\lambda^2(x) \geq \gamma) + \frac{V\gamma}{n\delta^2} \log^2 \frac{n}{\delta} + \frac{t}{n} \right).$$

REMARK. The following simple observation might be useful. Since

$$\mathbb{P}_n(\sigma_\lambda^2(x) \geq \gamma) \leq \frac{\mathbb{P}_n \sigma_\lambda^2}{\gamma},$$

one can plug this into the right-hand side of the bound of the theorem and then optimize it with respect to $\gamma$. The optimal value of $\gamma$ is

$$\hat{\gamma} := \frac{(\mathbb{P}_n \sigma_\lambda^2)^{1/2} \sqrt{n} \delta}{\sqrt{V} \log(n/\delta)} \wedge 1$$

(we are assuming here $\mathbb{P}_n \sigma_\lambda^2 > 0$!), which immediately leads to the following upper bound on generalization error:

$$K \inf_{0<\delta\leq\hat{\gamma}} \left( \mathbb{P}_n(yf_\lambda(x) \leq \delta) + 2\frac{\sqrt{V}(\mathbb{P}_n\sigma_\lambda^2)^{1/2}}{\sqrt{n}\delta} \log \frac{n}{\delta} \wedge \frac{V}{n\delta^2} \log^2 \frac{n}{\delta} + \frac{t}{n} \right).$$

This is to be compared with the bound (2.1) and it shows that the quantity $\mathbb{P}_n \sigma_\lambda^2$ might provide an interesting choice of complexity penalty in classification problems of this type. More generally, for $p \geq 1$ and (again, under the assumption $\mathbb{P}_n \sigma_\lambda^{2p} > 0$)

$$\hat{\gamma} := \frac{(\mathbb{P}_n\sigma_\lambda^{2p})^{1/(p+1)} n^{1/(p+1)} \delta^{2/(p+1)}}{V^{1/(p+1)} \log^{2/(p+1)}(n/\delta)} \wedge 1,$$

we are getting the bound

$$K \inf_{0<\delta\leq\hat{\gamma}} \left( \mathbb{P}_n(yf_\lambda(x) \leq \delta) \right.$$
$$\left. + 2\frac{V^{p/(p+1)}(\mathbb{P}_n\sigma_\lambda^{2p})^{1/(p+1)}}{n^{p/(p+1)} \delta^{2p/(p+1)}} \log^{2p/(p+1)} \frac{n}{\delta} \wedge \frac{V}{n\delta^2} \log^2 \frac{n}{\delta} + \frac{t}{n} \right).$$

In the limit $p \to \infty$ this yields the bound [provided that $\max_{1\leq j\leq n} \sigma_\lambda^2(X_j) > 0$]

$$K \inf_{0<\delta\leq\max_{1\leq j\leq n}\sigma_\lambda^2(X_j)} \left( \mathbb{P}_n(yf_\lambda(x) \leq \delta) + \frac{V\max_{1\leq j\leq n}\sigma_\lambda^2(X_j)}{n\delta^2} \log^2 \frac{n}{\delta} + \frac{t}{n} \right)$$

[which should be compared with (2.11); note the presence of the variance in the numerator].

The result of Theorem 3 is, probably, of limited interest since there is no reason to expect that the "global variances" of convex combinations output by popular learning algorithms are necessarily small. It is much more likely that it would be possible to split a convex combination into several clusters, each having a small variance. This is reflected in the following definition.

Given $m \geq 1$ and $\lambda \in \mathcal{P}(\mathcal{H})$, define a set

$$\mathcal{C}^m(\lambda) = \left\{ (\alpha_1, \ldots, \alpha_m, \lambda^1, \ldots, \lambda^m) : \lambda^k \in \mathcal{P}(\mathcal{H}), \alpha_k \geq 0, \sum_{k=1}^m \alpha_k \lambda^k = \lambda \right\}.$$

For an element $c \in \mathcal{C}^m(\lambda)$, we define a weighted variance over clusters by

$$(2.13) \qquad \sigma^2(c;x) = \sum_{k=1}^m \alpha_k^2 \sigma_{\lambda^k}^2(x),$$

where $\sigma_{\lambda^k}^2(x)$ are defined in (2.12). If indeed there are $m$ small clusters among functions $h_1, \ldots, h_T$, then one should be able to choose an element $c \in \mathcal{C}^m(\lambda)$ so that $\sigma^2(c;x)$ will be small on the majority of data points $X_1, \ldots, X_n$.

THEOREM 4. *If* (2.2) *holds, then there exists an absolute constant* $K > 0$ *such that for all* $t > 0$ *with probability at least* $1 - e^{-t}$ *for all* $\lambda \in \mathcal{P}(\mathcal{H})$ *and* $f(x) = f_\lambda(x) = \int h(x) \lambda(dh)$,

$$\mathbb{P}(y f_\lambda(x) \leq 0)$$
$$\leq K \inf_{m \geq 1} \inf_{c \in \mathcal{C}^m(\lambda)} \inf_{0 < \delta \leq \gamma \leq 1} \left( \mathbb{P}_n(y f_\lambda(x) \leq \delta) \right.$$
$$\left. + \mathbb{P}_n(\sigma^2(c; x) \geq \gamma) + \frac{V m \gamma}{n \delta^2} \log^2 \frac{n}{\delta} + \frac{t}{n} \right).$$

If we define the number of $(\gamma, \delta)$-clusters of $\lambda$ as the smallest $m$ for which there exists $c \in \mathcal{C}_\lambda$ such that

$$\mathbb{P}_n(\sigma^2(c; x) \geq \gamma) \leq \frac{V m \gamma}{n \delta^2} \log^2 \frac{n}{\delta}$$

and denote this number by $\hat{m}_\lambda(n, \gamma, \delta)$, then the bound implies that for all $\lambda \in \mathcal{P}(\mathcal{H})$

$$\mathbb{P}(y f_\lambda(x) \leq 0) \leq K \inf_{0 < \delta \leq \gamma} \left( \mathbb{P}_n(y f_\lambda(x) \leq \delta) + \frac{V \hat{m}_\lambda(n, \gamma, \delta) \gamma}{n \delta^2} \log^2 \frac{n}{\delta} + \frac{t}{n} \right).$$

The choice of $\gamma = \delta$ gives an upper bound with the error term (added to the empirical margin distribution) of the order

$$\frac{\hat{m}_\lambda(n, \delta, \delta)}{n \delta} \log^2 \frac{n}{\delta},$$

which significantly improves earlier bounds provided that we are lucky to have a small number of clusters $\hat{m}_\lambda(n, \delta, \delta)$ in the convex combination.

We now turn to a different approach to measuring complexity of convex combinations. It is based on empirical covering numbers of the set of functions involved in a particular convex combination. Let $\mathcal{H}$ be a class of measurable functions (classifiers) from $\mathcal{X}$ into $\{-1, 1\}$, such that $\mathcal{H}$ satisfies (2.2). It is interesting to note that in this case the condition (2.2) is equivalent to the condition that the class of sets $\mathcal{C} := \{\{h = +1\} : h \in \mathcal{H}\}$ is Vapnik–Chervonenkis (see, e.g., [13]).

As before, $\mathcal{H}$ will play the role of a base class. Let $\mathcal{F} := \mathrm{sconv}(\mathcal{H})$, that is, $\mathcal{F}$ is the symmetric convex hull of $\mathcal{H}$,

$$\mathrm{sconv}(\mathcal{H}) := \left\{ \sum_{i=1}^N \lambda_i h_i, h_i \in \mathcal{H}, \lambda_i \in \mathbb{R}, \sum_{i=1}^N |\lambda_i| \leq 1, N \geq 1 \right\}.$$

For $f \in \mathcal{F}$, a probability measure $Q$ on $\mathcal{X}$ and $p \in [1, +\infty]$, define

$$N_{d_{Q,p}}(f, \varepsilon) := \inf\{N_{d_{Q,p}}(\mathcal{H}', \varepsilon) : \mathcal{H}' \subset \mathcal{H}, f \in \mathrm{sconv}(\mathcal{H}')\}.$$

Let us call a subset $\mathcal{H}' \subset \mathcal{H}$ *a base* of $f \in \operatorname{sconv}(\mathcal{H})$ iff $f \in \operatorname{sconv}(\mathcal{H}')$. Then $N_{d_{Q,p}}(f; \varepsilon)$ is the minimal $\varepsilon$-covering number of bases of $f$. Let

$$\hat{\psi}_n(f; \delta) := \int_0^\delta \sqrt{N_{d_{\mathbb{P}_n, 2}}(f, \varepsilon) \log(1/\varepsilon)}\, d\varepsilon.$$

As earlier in this section (see also [21]), for a concave nondecreasing function $\psi$ on $[0, +\infty)$ with $\psi(0) = 0$, we define $\varepsilon_n^\psi(\delta)$ as the largest solution of the equation

$$\varepsilon = \frac{1}{\delta \sqrt{n}} \psi(\delta \sqrt{\varepsilon})$$

with respect to $\varepsilon$. Let now

$$\hat{\varepsilon}_n(f, \delta) := \varepsilon_n^{\hat{\psi}_n(f, \cdot)}(\delta).$$

The function $\hat{\psi}_n(f, \cdot)$ can be viewed as a data- and classifier-dependent estimate of the entropy integral in the condition (2.6), and the bound of Theorem 5 below is an adaptive version of $\psi$-bounds developed in [21].

THEOREM 5. *If a class of measurable functions $\mathcal{H} = \{h \colon \mathcal{X} \to \{-1, +1\}\}$ satisfies* (2.2), *then for all $t \geq C \log^2 n$, with probability at least $1 - e^{-t}$ the following bound holds for all $f \in \mathcal{F}$:*

$$\mathbb{P}\{yf(x) \leq 0\} \leq K \inf_{\delta \in (0,1]} \left[ \mathbb{P}_n\{yf(x) \leq \delta\} + \hat{\varepsilon}_n(f, \delta) + \frac{t}{n\delta^2} \right],$$

*where $K, C > 0$ are absolute constants.*

REMARK 1. Clearly, for all $\varepsilon > 0$

$$N_{d_{\mathbb{P}_n, 2}}(f, \varepsilon) \leq N_{d_{\mathbb{P}_n, \infty}}(f, \varepsilon),$$

and since the functions in $\mathcal{H}$ take their values in $\{-1, 1\}$, $N_{d_{\mathbb{P}_n, \infty}}(f, \varepsilon)$ does not depend on $\varepsilon$ for all $\varepsilon < 2$. Therefore, in this range of $\varepsilon$ we will use the notation $N_{d_{\mathbb{P}_n, \infty}}(f)$ for it. This quantity is always bounded by $2^n$ and it shows how many classifiers $h_j \in \mathcal{H}$ that differ on the sample are involved in the "most economical" representation of $f \in \operatorname{sconv}(\mathcal{H})$ (so it can be viewed as a dimension of $f$). The following bound is trivial:

$$\hat{\psi}_n(f; \delta) \leq 2\sqrt{N_{d_{\mathbb{P}_n, \infty}}(f)} \delta \sqrt{\log \frac{1}{\delta}}, \qquad \delta < e^{-1},$$

and it shows, in particular, that $\hat{\psi}_n(f, \delta)$ is well defined. It also shows that the function $\hat{\varepsilon}_n(f, \delta)$ involved in the bound of the theorem can be replaced by the following upper bound that has a much simpler meaning:

$$\frac{N_{d_{\mathbb{P}_n, \infty}}(f)}{n} \log \frac{n}{\delta N_{d_{\mathbb{P}_n, \infty}}(f)}$$

[although $\hat{\varepsilon}_n(f, \delta)$ can be much smaller than this upper bound].

REMARK 2. In fact, the bound of the theorem can be improved by introducing

$$\hat{H}_n(f,\varepsilon) := N_{d_{\mathbb{P}_n,2}}(f,\varepsilon)\log\frac{1}{\varepsilon} \wedge \varepsilon^{-2V/(V+2)}$$

and defining the function

$$\hat{\psi}_n(f,t,\delta) := \int_0^\delta \hat{H}_n^{1/2}\left(f,\varepsilon \vee \sqrt{\frac{t}{n}}\right) d\varepsilon.$$

Then one can define $\hat{\varepsilon}_n(f,t,\delta)$ as $\varepsilon_n^\psi(\delta)$ with $\psi(\cdot) := \hat{\psi}_n(f,t,\cdot)$. It follows from the proofs below that for all $t \geq C\log^2 n$, with probability at least $1 - e^{-t}$ the following bound holds for all $f \in \mathcal{F}$:

$$\mathbb{P}\{yf(x) \leq 0\} \leq K \inf_{\delta \in (0,1]}\left[\mathbb{P}_n\{yf(x) \leq \delta\} + \hat{\varepsilon}_n(f,t,\delta) + \frac{t}{n\delta^2}\right]$$

with some constants $K, C > 0$. The term $\varepsilon^{-2V/(V+2)}$ in the definition of $\hat{H}_n(f,\varepsilon)$ is (up to a constant) a well-known upper bound on the entropy of the convex hull of a VC-type class. The definition of $\hat{H}_n(f,\varepsilon)$ is based on an upper bound (see Lemma 2 below) on the entropy of the *restricted convex hull* of $\mathcal{H}$ defined (given a probability measure $Q$ and $p \geq 1$) as

$$\{f \in \mathrm{sconv}(\mathcal{H}) : \forall \varepsilon : N_{d_{Q,p}}(f,\varepsilon) \leq N(\varepsilon)\},$$

where $N$ is a given nonincreasing function. In fact, any other upper bound on the entropy of such sets can be used instead of $\hat{H}_n(f,\varepsilon)$. Apparently, more subtle bounds than the result of Lemma 2 (that interpolate better between the case of finite-dimensional convex combinations and the case of the whole convex hull) should exist and allow one to improve the bound of Theorem 5, but at the moment we do not know how to prove a better bound. Theorem 5 can be extended to classes $\mathcal{H}$ of functions taking values in $[-1, 1]$ (not necessarily binary functions), but its formulation becomes more complicated since it involves both $\mathcal{L}_2(\mathbb{P}_n)$- and $\mathcal{L}_1(\mathbb{P}_n)$-entropies in this case.

**3. Proofs.** Theorem 6 will be the main technical tool in the proofs of Theorems 1–4. This theorem extends the inequality of Vapnik and Chervonenkis for VC-classes of sets and VC-major classes of functions to classes of functions $\mathcal{F} = \{f : \mathcal{X} \to [-1,1]\}$ satisfying the uniform entropy condition

(3.1) $$\int_0^\infty \log^{1/2} N(\mathcal{F},u)\,du < \infty,$$

where

$$N(\mathcal{F},u) = \sup_{Q \in \mathcal{P}(\mathcal{X})} N_{d_{Q,2}}(\mathcal{F},u).$$

For instance, it obviously holds under (2.2) for $\mathcal{F} = \mathrm{conv}(\mathcal{H})$, as it follows from the well-known bounds on the entropy of the convex hull.

THEOREM 6. *If $\mathcal{F} = \{f : \mathcal{X} \to [0,1]\}$ is a class of $[0,1]$-valued functions that satisfies (3.1), then there exists an absolute constant $K > 0$ such that for any $t > 0$ with probability at least $1 - e^{-t}$ for all $f \in \mathcal{F}$*

$$(3.2) \quad \mathbb{P}f - \mathbb{P}_n f \leq K \left( n^{-1/2} \int_0^{(\mathbb{P}f)^{1/2}} \log^{1/2} N(\mathcal{F}, u) \, du + \left( \frac{t\mathbb{P}f}{n} \right)^{1/2} \right),$$

*and with probability at least $1 - e^{-t}$ for all $f \in \mathcal{F}$*

$$(3.3) \quad \mathbb{P}_n f - \mathbb{P}f \leq K \left( n^{-1/2} \int_0^{(\mathbb{P}_n f)^{1/2}} \log^{1/2} N(\mathcal{F}, u) \, du + \left( \frac{t\mathbb{P}_n f}{n} \right)^{1/2} \right).$$

PROOF. Equation (3.2) is Corollary 1 in [25]. Equation (3.3) is not formulated in [25] explicitly but it is proved similarly to (3.2). Equations (3.2) and (3.3) also follow easily from Corollary 3 in [26]. □

There are two features of this result that make it particularly useful. First of all, it is well known (see [13]) that if, given $p > 0$, we look at the layer of functions $\mathcal{F}_p = \{f \in \mathcal{F} : \mathbb{P}f \leq p\}$, then the typical value of the deviation $\mathbb{P}f - \mathbb{P}_n f$ on this layer or, in other words, the expectation $\mathbb{E}\sup\{\mathbb{P}f - \mathbb{P}_n f : f \in \mathcal{F}_p\}$, can be estimated by the entropy integral

$$n^{-1/2} \int_0^{\sqrt{p}} \log^{1/2} N(\mathcal{F}, u) \, du,$$

where the upper limit $\sqrt{p}$ measures the size of $\mathcal{F}_p$. This simply reflects the fact that functions with smaller mean $\mathbb{P}f$ will have smaller fluctuations. Theorem 6 says that this happens on all layers at the same time, which gives us an adaptive control over the whole class $\mathcal{F}$. The second important feature of this result is that the deviation from a typical value is controlled for each function individually by the term $(t\mathbb{P}f/n)^{1/2}$. This is convenient from the point of view of structural risk minimization since one only has to estimate the typical value on each class to which a function $f$ may belong, but the deviation term is left unchanged. For other results in this direction we refer the reader to [26].

Given an integer $d \geq 1$, denote

$$\mathcal{F}_d = \mathrm{conv}_d(\mathcal{H}) = \left\{ \sum_{i=1}^d \lambda_i h_i : \sum_{i=1}^d \lambda_i \leq 1, \lambda_i \geq 0, h_i \in \mathcal{H} \right\}.$$

Again, let $\Phi = \{\varphi_\delta : \mathbb{R} \to [0,1] : \delta \in \Delta \subset \mathbb{R}_+\}$ be a countable family of Lipschitz functions such that Lipschitz norm of $\varphi_\delta$ is equal to $\delta^{-1}$ and $\sum_{\delta \in \Delta} \delta < \infty$. One can use a specific choice of $\Delta = \{2^{-k} : k \geq 1\}$. For $a > 0, b \geq 0$ we define

$$\phi(a, b) = \frac{(a-b)^2}{a} I(a \geq b),$$

and for $a = 0$ we let $\phi(a, b) = \phi(0, b) = 0$. The following theorem holds.

THEOREM 7. *If* (2.2) *holds, then there exists* $K > 0$ *such that for all* $t > 0$ *with probability at least* $1 - e^{-t}$ *we have for all* $d \geq 1, f \in \mathcal{F}_d$ *and* $\delta \in \Delta$,

$$(3.4) \qquad \phi(\mathbb{P}\varphi_\delta(yf(x)), \mathbb{P}_n\varphi_\delta(yf(x))) \leq K\left(\frac{dV}{n}\log\frac{n}{\delta} + \frac{t}{n}\right).$$

PROOF. The proof is a straightforward application of Theorem 6. We will proceed in several steps.

*Step* 1 (Estimating covering numbers). First of all, if given a class of measurable functions on $\mathcal{X}$, $\mathcal{F} = \{f : \mathcal{X} \to [0,1]\}$, we introduce a new class of measurable functions

$$\mathcal{F}^{\mathcal{Y}} = \{g(x,y) = yf(x) : \mathcal{X} \times \mathcal{Y} \to [-1,1] : f \in \mathcal{F}\}$$

defined on $\mathcal{X} \times \mathcal{Y}$, then

$$N(\mathcal{F}^{\mathcal{Y}}, u) = N(\mathcal{F}, u)$$

since for any $(x_1, y_1), \ldots, (x_n, y_n)$ and any $f_1, f_2 \in \mathcal{F}$ we have

$$\frac{1}{n}\sum_{i=1}^{n}(y_i f_1(x_i) - y_i f_2(x_i))^2 = \frac{1}{n}\sum_{i=1}^{n}(f_1(x_i) - f_2(x_i))^2.$$

Therefore, condition (2.2) on $\mathcal{H}$ is equivalent to the corresponding condition on $\mathcal{H}^{\mathcal{Y}}$.

The following bound for the uniform entropy of $\mathcal{F}_d^{\mathcal{Y}}$ in terms of $N(\mathcal{H}^{\mathcal{Y}}, u)$ is well known (see [21], Lemma 2):

$$N(\mathcal{F}_d^{\mathcal{Y}}, u) \leq \left(\frac{2e^2 N(\mathcal{H}^{\mathcal{Y}}, u)(d^2 + 16u^{-2})}{d^2}\right)^d.$$

In combination with (2.2) it implies that for some $K > 0$

$$\log N(\mathcal{F}_d^{\mathcal{Y}}, u) \leq KdV \log\frac{1}{u}.$$

For a fixed $\varphi_\delta \in \Phi$ the uniform covering numbers of the class $\varphi_\delta \circ \mathcal{F}_d^{\mathcal{Y}} = \{\varphi_\delta(g) : g \in \mathcal{F}_d^{\mathcal{Y}}\}$ can be bounded as

$$N(\varphi_\delta \circ \mathcal{F}_d^{\mathcal{Y}}, u) \leq N(\mathcal{F}_d^{\mathcal{Y}}, \delta u),$$

since for any probability measure $Q$ on $\mathcal{X} \times \mathcal{Y}$ the Lipschitz condition on $\varphi_\delta$ implies that

$$(Q(\varphi_\delta(yf(x)) - \varphi_\delta(yg(x)))^2)^{1/2} \leq \delta^{-1}(Q(f-g)^2)^{1/2},$$

and, therefore,

$$\log N(\varphi_\delta \circ \mathcal{F}_d^{\mathcal{Y}}, u) \leq KdV \log\frac{1}{\delta u}.$$

*Step* 2 (Nonadaptive bound). Theorem 6 applied to $\varphi_\delta \circ \mathcal{F}_d^{\mathcal{Y}}$ guarantees that for any $t > 0$ with probability at least $1 - e^{-t}$ for all $f \in \mathcal{F}_d$,

$$\mathbb{P}\varphi_\delta(yf(x)) - \mathbb{P}_n\varphi_\delta(yf(x))$$
$$\leq K\left(\left(\frac{dV}{n}\right)^{1/2}\int_0^{(\mathbb{P}\varphi_\delta(yf(x)))^{1/2}} \log^{1/2}\frac{1}{\delta u}\,du + \left(\frac{t\mathbb{P}\varphi_\delta(yf(x))}{n}\right)^{1/2}\right).$$

To estimate the first term on the right-hand side one can easily check that

$$(3.5) \qquad \int_0^s \left(\log\frac{1}{u}\right)^{1/2}du \leq 2s\left(\log\frac{1}{s}\right)^{1/2} \qquad \text{for } s \in [0, e^{-1}].$$

This inequality is well known and, moreover, the value 2 of the constant is irrelevant here. Hence,

$$\int_0^{(\mathbb{P}\varphi_\delta(yf(x)))^{1/2}} \log^{1/2}\frac{1}{\delta u}\,du$$
$$= \delta^{-1}\int_0^{\delta(\mathbb{P}\varphi_\delta(yf(x)))^{1/2}} \log^{1/2}\frac{1}{s}\,ds$$
$$\leq 2(\mathbb{P}\varphi_\delta(yf(x)))^{1/2}\max\left(1, \log^{1/2}\frac{1}{\delta(\mathbb{P}\varphi_\delta(yf(x)))^{1/2}}\right).$$

Without loss of generality we can assume that $\mathbb{P}\varphi_\delta(yf(x)) \geq n^{-1}$; otherwise, the bound of the theorem becomes trivial. Therefore,

$$\max\left(1, \log^{1/2}\frac{1}{\delta(\mathbb{P}\varphi_\delta(yf(x)))^{1/2}}\right) \leq \log^{1/2}\frac{n}{\delta},$$

which finally yields

$$\int_0^{(\mathbb{P}\varphi_\delta(yf(x)))^{1/2}} \log^{1/2}\frac{1}{\delta u}\,du \leq 2(\mathbb{P}\varphi_\delta(yf(x)))^{1/2}\log^{1/2}\frac{n}{\delta}.$$

We have proved that

$$\frac{\mathbb{P}\varphi_\delta(yf(x)) - \mathbb{P}_n\varphi_\delta(yf(x))}{(\mathbb{P}\varphi_\delta(yf(x)))^{1/2}} \leq K\left(\left(\frac{dV}{n}\log\frac{n}{\delta}\right)^{1/2} + \left(\frac{t}{n}\right)^{1/2}\right),$$

which implies that

$$\phi(\mathbb{P}\varphi_\delta(yf(x)), \mathbb{P}_n\varphi_\delta(yf(x))) \leq K\left(\frac{dV}{n}\log\frac{n}{\delta} + \frac{t}{n}\right).$$

*Step* 3 (Union bound, adaptivity). The statement of the theorem now follows by applying the union bound and increasing $K$. Indeed, let us introduce the event

$$A_{d,\delta}(t') = \left\{\forall f \in \mathcal{F}_d : \phi(\mathbb{P}\varphi_\delta(yf(x)), \mathbb{P}_n\varphi_\delta(yf(x))) \leq K\left(\frac{dV}{n}\log\frac{n}{\delta} + \frac{t'}{n}\right)\right\},$$

which holds with probability $1 - e^{-t'}$. For a fixed $t$ and for a fixed $d$ and $\delta$, define $t'$ according to the equality $e^{-t'} = (\delta e^{-t})/(d^2 K)$, where $K$ is chosen so that the condition $\sum_{d \in Z_+, \delta \in \Delta} \delta d^{-2}/K \leq 1$ holds. With this choice of $t'$ the event $A_{d,\delta}(t')$ can be rewritten

$$A_{d,\delta} = \Big\{ \forall f \in \mathcal{F}_d : \phi(\mathbb{P}\varphi_\delta(yf(x)), \mathbb{P}_n \varphi_\delta(yf(x)))$$

$$\leq K\Big(\frac{dV}{n}\log\frac{n}{\delta} + \frac{1}{n}\log\frac{Kd^2}{\delta} + \frac{t}{n}\Big)\Big\},$$

and its probability is greater than

$$\mathbf{Pr}(A_{d,\delta}) \geq 1 - \frac{\delta e^{-t}}{d^2 K}.$$

It implies that the probability of the intersection

$$\mathbf{Pr}\Big(\bigcap_{d,\delta} A_{d,\delta}\Big) \geq 1 - \sum_{\delta,d} \frac{\delta e^{-t}}{d^2 K} \geq 1 - e^{-t}.$$

This means that with probability at least $1 - e^{-t}$ all the events $A_{d,\delta}$ hold simultaneously. But, obviously, the second term in the definition of $A_{d,\delta}$ can be bounded by

$$\frac{1}{n}\log\frac{Kd^2}{\delta} \leq K\frac{d}{n}\log\frac{n}{\delta}$$

and, thus, $A_{\delta,d}$ is a subset of the event

$$A_{d,\delta} \subseteq A'_{d,\delta} = \Big\{\forall f \in \mathcal{F}_d : \phi(\mathbb{P}\varphi_\delta(yf(x)), \mathbb{P}_n \varphi_\delta(yf(x))) \leq K'\Big(\frac{dV}{n}\log\frac{n}{\delta} + \frac{t}{n}\Big)\Big\}$$

for some $K' > K$, which proves the statement of the theorem, since

$$\mathbf{Pr}\Big(\bigcap_{d,\delta} A'_{d,\delta}\Big) \geq \mathbf{Pr}\Big(\bigcap_{d,\delta} A_{d,\delta}\Big) \geq 1 - e^{-t}.$$

□

PROOF OF THEOREM 1. For a fixed $d, \gamma$ consider a class $\mathcal{F}_{d,\gamma} = \{f \in \mathcal{F} : d(f;\gamma) \leq d\}$. One can estimate the uniform entropy of $\mathcal{F}_{d,\gamma}$ as (see [21])

$$\log N(\mathcal{F}_{d,\gamma}, u) \leq K\Big(d\log\frac{1}{u} + \Big(\frac{\gamma}{u}\Big)^\alpha\Big).$$

For a fixed $\varphi_\delta \in \Phi$ the uniform covering numbers of the class $\varphi_\delta \circ \mathcal{F}^{\mathcal{Y}}_{d,\gamma} = \{\varphi_\delta(yf(x)) : f \in \mathcal{F}_{d,\gamma}\}$ can be bounded as

$$N(\varphi_\delta \circ \mathcal{F}^{\mathcal{Y}}_{d,\gamma}, u) \leq N(\mathcal{F}_{d,\gamma}, \delta u),$$

since, for any probability measure $Q$ on $\mathcal{X} \times \mathcal{Y}$, the Lipschitz condition on $\varphi_\delta$ implies that

$$(Q(\varphi_\delta(yf(x)) - \varphi_\delta(yg(x)))^2)^{1/2} \leq \delta^{-1}(Q(yf - yg)^2)^{1/2} = \delta^{-1}(Q(f-g)^2)^{1/2},$$

and, therefore,

$$\log N(\varphi_\delta \circ \mathcal{F}_{d,\gamma}^{\mathcal{Y}}, u) \leq K\left(d\log\frac{1}{\delta u} + \left(\frac{\gamma}{\delta u}\right)^\alpha\right).$$

Using this estimate on the covering numbers, Theorem 6 now implies (in exactly the same way we used it in the proof of Theorem 7; only integration here is easier) that for any $t > 0$ with probability at least $1 - e^{-t}$ for all $f \in \mathcal{F}_{d,\gamma}$

$$\frac{\mathbb{P}\varphi_\delta(yf(x)) - \mathbb{P}_n\varphi_\delta(yf(x))}{(\mathbb{P}\varphi_\delta(yf(x)))^{1/2}}$$

$$\leq K\left(\left(\frac{d}{n}\log\frac{n}{\delta}\right)^{1/2} + \left(\frac{\gamma}{\delta}\right)^{\alpha/2}\frac{(\mathbb{P}\varphi_\delta(yf(x)))^{-\alpha/4}}{n^{1/2}} + \left(\frac{t}{n}\right)^{1/2}\right).$$

It remains to show that, possibly increasing $K$, this inequality holds for all $d, \delta$ and $\gamma$. To do this we will use the above inequality with $t$ replaced by $t' + \log\frac{Kd^2}{\delta\gamma}$ and, hence, $e^{-t}$ replaced by $e^{-t'} = (e^{-t}\delta\gamma)/(Kd^2)$, where $\delta, \gamma \in \{2^{-k} : k \geq 1\}$. Then the union bound should be applied in the whole range of $d, \delta$ and $\gamma$. Without loss of generality we assume that for all $f \in \mathcal{F}$ and $\delta \in \Delta$ we have $\mathbb{P}\varphi_\delta(yf(x)) \geq n^{-1}$, and $\gamma$ can be restricted to the set of values satisfying

$$\left(\frac{\gamma}{\delta}\right)^{\alpha/2}\frac{(\mathbb{P}\varphi_\delta(yf(x)))^{-\alpha/4}}{n^{1/2}} \geq \left(\frac{t}{n}\right)^{1/2},$$

or, equivalently,

$$\gamma \geq \delta(\mathbb{P}\varphi_\delta(yf(x)))^{1/2}t^{1/\alpha} \geq \delta n^{-1/2}t^{1/\alpha}.$$

Under these assumptions

$$\log\frac{Kd^2}{\delta\gamma} \leq Kd\log\frac{n}{\delta},$$

which allows us to complete the proof by using the union bound and choosing the value of $K$ large enough. $\square$

PROOF OF COROLLARY 1. To see that Theorem 1 implies Corollary 1 one should first notice that if $\mathbb{P}_n\varphi_\delta(yf(x)) = 0$, then the inequality of Theorem 1 can be solved for $\mathbb{P}\varphi_\delta(yf(x))$ to give

$$\mathbb{P}\varphi_\delta(yf(x)) \leq I(f) = K\inf_\gamma\left(\frac{d(f;\gamma)}{n}\log\frac{n}{\delta} + \left(\frac{\gamma}{\delta}\right)^{2\alpha/(2+\alpha)}n^{-2/(\alpha+2)} + \frac{t}{n}\right)$$

(we prove it below). Moreover, if $\mathbb{P}_n \varphi_\delta(yf(x))$ is of the same order of magnitude as $I(f)$, then we will show that $\mathbb{P}\varphi_\delta(yf(x))$ will also be of the same order of magnitude as $I(f)$. Finally, if $\mathbb{P}_n \varphi_\delta(yf(x))$ is larger than a constant times $I(f)$, then $\mathbb{P}\varphi_\delta(yf(x))$ is dominated by a constant times $\mathbb{P}_n \varphi_\delta(yf(x))$. After all this is proved, it remains to notice that, for a specific choice of functions $\varphi_\delta$ such that $\varphi_\delta(s) = 1$ for $s \leq 0$, $\varphi_\delta(s) = 0$ for $s \geq \delta$ and linear on $[0, \delta]$, we have

$$\mathbb{P}(yf(x) \leq 0) \leq \mathbb{P}\varphi_\delta(yf(x)) \quad \text{and} \quad \mathbb{P}_n\varphi_\delta(yf(x)) \leq \mathbb{P}_n(yf(x) \leq \delta).$$

We will now explain how to solve the inequality of Theorem 1. We observe that it is of the form

$$(3.6) \qquad y \leq x + ay^{1/2} + by^\beta,$$

where $y = \mathbb{P}\varphi_\delta, x = \mathbb{P}_n \varphi_\delta$, $0 < \beta < 1, a, b > 0$. In our case also $\beta = 1/2 - \alpha/4$. Define $y_1$ and $y_2$ as the solutions of the equations

$$y_1 = ay_1^{1/2}, \qquad y_2 = by_2^\beta$$

and notice that

$$y \geq ay^{1/2} \quad \text{for } y \geq y_1; \qquad y \geq by^\beta \quad \text{for } y \geq y_2.$$

Assume that $x \leq y_1 + y_2$. Then (3.6) implies that $y \leq K(y_1 + y_2)$ for some absolute constant $K > 0$. Indeed, if we plug $K(y_1 + y_2)$ into the right-hand side of (3.6) we get

$$x + a(K(y_1 + y_2))^{1/2} + b(K(y_1 + y_2))^\beta$$
$$\leq (y_1 + y_2) + K^{1/2}a(y_1 + y_2)^{1/2} + K^\beta b(y_1 + y_2)^\beta$$
$$\leq (y_1 + y_2) + K^{1/2}(y_1 + y_2) + K^\beta(y_1 + y_2)^\beta$$
$$\qquad \text{(since } y_1 + y_2 \geq y_1 \text{ and } y_1 + y_2 \geq y_2\text{)}$$
$$\leq (1 + K^{1/2} + K^\beta)(y_1 + y_2) \leq K(y_1 + y_2),$$

if $K$ is large enough. This shows that (3.6) fails for $y \geq K(y_1 + y_2)$, and hence the solution of (3.6) is smaller than $K(y_1 + y_2)$. Assuming that $x \geq y_1 + y_2$ and setting $C := \frac{y}{x}$, we get from (3.6)

$$Cx \leq x + C^{1/2}ax^{1/2} + C^\beta bx^\beta \leq x + C^{1/2}x + C^\beta x = (1 + C^{1/2} + C^\beta)x,$$

which implies $C \leq 1 + C^{1/2} + C^\beta$ and hence $y \leq Kx$ for a large enough constant $K$. Thus, always with large enough $K$ we have $y \leq K(x + y_1 + y_2)$, implying the result. $\square$

PROOF OF THEOREM 2. Let us make a specific choice of functions $\varphi_\delta$. For each $\delta \in \Delta$ we set $\varphi_\delta$ to be $\varphi_\delta(s) = 1$ for $s \leq \delta$, $\varphi_\delta(s) = 0$ for $s \geq 2\delta$ and linear on $[\delta, 2\delta]$.

Let us fix $f = \sum_{k=1}^{T} \lambda_k h_k \in \mathcal{F}$, and for a fixed $0 \leq d \leq T$ represent $f$ as

$$f = \sum_{k=1}^{d} \lambda_k h_k + \gamma_d(f) \sum_{k=d+1}^{T} \lambda'_k h_k,$$

where $\gamma_d(f) = \sum_{k=d+1}^{T} \lambda_k$ and $\lambda'_k = \lambda_k/\gamma_d(f)$.

Given $N \geq 1$, we generate an i.i.d. sequence of functions $\xi_1, \ldots, \xi_N$ according to the distribution $\mathbb{P}_\xi(\xi_i = h_k) = \lambda'_k$ for $k = d+1, \ldots, T$ and independent of $\{(X_k, Y_k)\}$. Clearly, $\mathbb{E}_\xi \xi_i(x) = \sum_{k=d+1}^{T} \lambda'_k h_k(x)$. Consider a function

$$g(x) = \sum_{k=1}^{d} \lambda_k h_k(x) + \gamma_d(f) \frac{1}{N} \sum_{k=1}^{N} \xi_k(x),$$

which plays the role of a random approximation of $f$ in the following sense. We can write

$$\mathbb{P}(yf(x) \leq 0) = \mathbb{E}_\xi \mathbb{P}(yf(x) \leq 0, yg(x) \leq \delta) + \mathbb{E}_\xi \mathbb{P}(yf(x) \leq 0, yg(x) \geq \delta)$$
(3.7)
$$\leq \mathbb{E}_\xi \mathbb{P} \varphi_\delta(yg(x)) + \mathbb{E}\mathbb{P}_\xi(yg(x) \geq \delta, \mathbb{E}_\xi yg(x) \leq 0).$$

In the last term for a fixed $(x,y) \in \mathcal{X} \times \mathcal{Y}$ we have

$$\mathbb{P}_\xi(yg(x) \geq \delta, \mathbb{E}_\xi yg(x) \leq 0) \leq \mathbb{P}_\xi(yg(x) - \mathbb{E}_\xi yg(x) \geq \delta)$$

$$= \mathbb{P}_\xi\left(\sum_{i=1}^{N}(y\xi_i(x) - y\mathbb{E}_\xi \xi_i(x)) \geq N\delta/\gamma_d(f)\right)$$

$$\leq \exp(-N\delta^2/2\gamma_d^2(f)),$$

where in the last step we used Hoeffding's inequality. Hence,

(3.8) $$\mathbb{P}(yf(x) \leq 0) - e^{-N\delta^2/2\gamma_d^2(f)} \leq \mathbb{E}_\xi \mathbb{P} \varphi_\delta(yg(x)).$$

Similarly, one can write

$$\mathbb{E}_\xi \mathbb{P}_n \varphi_\delta(yg(x)) \leq \mathbb{E}_\xi \mathbb{P}_n(yg(x) \leq 2\delta) \leq \mathbb{P}_n(yf(x) \leq 3\delta)$$
(3.9)
$$+ \mathbb{E}_\xi \mathbb{P}_n(yg(x) \leq 2\delta, yf(x) \geq 3\delta)$$

$$\leq \mathbb{P}_n(yf(x) \leq 3\delta) + e^{-N\delta^2/2\gamma_d^2(f)}.$$

Clearly, for any random realization of the sequence $\xi_1, \ldots, \xi_N$, the random function $g$ belongs to the class $\mathcal{F}_{d+N}$. Convexity of the function $\phi(a,b)$ and Theorem 7 imply that for any $t > 0$ with probability at least $1 - e^{-t}$ for all $\delta \in \Delta$ and all $f \in \mathcal{F}$

$$\phi(\mathbb{E}_\xi \mathbb{P} \varphi_\delta(yg(x)), \mathbb{E}_\xi \mathbb{P}_n \varphi_\delta(yg(x))) \leq \mathbb{E}_\xi \phi(\mathbb{P} \varphi_\delta(yg(x)), \mathbb{P}_n \varphi_\delta(yg(x)))$$

$$\leq K\left(\frac{V(d+N)}{n} \log \frac{n}{\delta} + \frac{t}{n}\right).$$

The fact that $\phi(a,b)$ is decreasing in $b$ and increasing in $a$ combined with (3.8) and (3.9) implies that

$$\phi(\mathbb{P}(yf(x) \leq 0) - e^{-N\delta^2/2\gamma_d^2(f)}, \mathbb{P}_n(yf(x) \leq 3\delta) + e^{-N\delta^2/2\gamma_d^2(f)})$$
$$\leq K\bigg(\frac{V(d+N)}{n}\log\frac{n}{\delta} + \frac{t}{n}\bigg).$$

Setting $N = 2(\gamma_d^2(f)/\delta^2)\log n$, we get

$$\phi(\mathbb{P}(yf(x) \leq 0) - 1/n, \mathbb{P}_n(yf(x) \leq 3\delta) + 1/n) \leq K\bigg(\frac{Ve_n(f,\delta,d)}{n}\log\frac{n}{\delta} + \frac{t}{n}\bigg),$$

where $e_n(f,\delta,d) = d + 2(\gamma_d^2(f)/\delta^2)\log n$. Solving the last inequality for $\mathbb{P}(yf(x) \leq 0)$ and changing the variable $3\delta \to \delta$ gives the bound (that holds with probability at least $1 - e^{-t}$)

(3.10) $\qquad \mathbb{P}(yf(x) \leq 0) \leq (W^{1/2} + (\mathbb{P}_n(yf(x) \leq \delta) + W)^{1/2})^2,$

where

$$W = W(f,n,d,\delta,t) := K\bigg(\frac{Ve_n(f,\delta,d)}{n}\log\frac{n}{\delta} + \frac{t}{n}\bigg).$$

It remains to make the bound uniform over $d$ and $\delta$, which is done using standard union bound techniques. More specifically, replace $t$ in the above bound by $t'(d,\delta) = t + 2\log(1/\delta) + 2\log d + c$, where $\delta \in \{2^{-k} : k \geq 1\}$ and

$$c := 2\log\bigg(\sum_{k=1}^{\infty} k^{-2}\bigg).$$

Then the union bound can be used to show that (3.10) [with $t$ replaced by $t'(d,\delta)$] holds for all $d$ and all $\delta \in \{2^{-k} : k \geq 1\}$ simultaneously with probability at least $1 - p$, where

$$p \leq e^{-t-c} \sum_{k=1,d=1}^{\infty} e^{-2\log k - 2\log d} = e^{-t-c}\bigg(\sum_{k=1}^{\infty} k^{-2}\bigg)^2 = e^{-t},$$

and, hence, we also have with probability at least $1 - e^{-t}$

$$\mathbb{P}(yf(x) \leq 0) \leq \inf_{\delta \in \{2^{-k} : k \geq 1\}} \inf_d (W^{1/2}(f,n,d,\delta,t'(d,\delta))$$
$$+ (\mathbb{P}_n(yf(x) \leq \delta) + W(f,n,d,\delta,t'(d,\delta)))^{1/2})^2.$$

Taking into account the monotonicity of the function $e_n(f,\delta,d)$ with respect to $\delta$ (and increasing the value of the constant $K$), it is now easy to extend the infimum over $\delta$ to all $\delta \in (0,1]$. Increasing the value of $K$ further allows one to rewrite the bound as

$$\mathbb{P}(yf(x) \leq 0) \leq \inf_{\delta \in (0,1]} (U^{1/2} + (\mathbb{P}_n(yf(x) \leq \delta) + U)^{1/2})^2$$

with $U$ defined in the formulation of the theorem, which completes the proof. □

Theorem 3 is a special case of Theorem 4; thus we will proceed by proving Theorem 4.

PROOF OF THEOREM 4. We will proceed to prove Theorem 4 in several steps.

*Step* 1 (Random approximation). Consider functions $\varphi_\delta$ the same as in the proof of Theorem 2. Let $\lambda \in \mathcal{P}(\mathcal{H})$ and $f(x) = \int h(x)\lambda(dh)$. Consider an element $c \in \mathcal{C}^m(\lambda)$, that is, $c = (\alpha_1, \ldots, \alpha_m, \lambda^1, \ldots, \lambda^m)$, such that $\lambda = \sum_{i=1}^m \alpha_j \lambda^j$ and $\lambda^j \in \mathcal{P}(\mathcal{H})$. We interpreted $c$ as a decomposition of $\lambda$ into $m$ clusters, or in other words, the decomposition of the set $\{h_i\}$ into $m$ clusters. This time we will generate functions from each cluster independently from each other (and, as before, independently of the data) and take their weighted sum to approximate $f(x)$. Given $N \geq 1$, let us generate independent random functions $\xi_k^j(x)$, $k \leq N, j \leq m$, where for each $j \leq m$, the $\xi_k^j$'s have the distribution

$$\mathbb{P}_{\xi^j}(\xi_k^j = h_i) = \lambda^j(\{h_i\}) = \lambda_i^j, \qquad i \leq T.$$

Consider a function

$$g(x) = \frac{1}{N}\sum_{j=1}^m \alpha_j \sum_{k=1}^N \xi_k^j(x) = \frac{1}{N}\sum_{k=1}^N g_k(x),$$

where $g_k(x) = \sum_{j=1}^m \alpha_j \xi_k^j(x)$. For a fixed $x \in \mathcal{X}$ and $k \leq N$, the variance of $g_k$ with respect to the distribution $\mathbb{P}_\xi = \mathbb{P}_{\xi^1} \times \cdots \times \mathbb{P}_{\xi^m}$ is

$$\mathrm{Var}_\xi(g_k(x)) = \sum_{j=1}^m \alpha_j^2 \mathrm{Var}_\xi(\xi_1^j(x)) = \sum_{j=1}^m \alpha_j^2 \sigma_{\lambda^j}^2(x) = \sigma^2(c;x).$$

The main difference from the proof of Theorem 2 is that in (3.7) we also introduce the condition on the variance $\sigma^2(c;x)$. Namely, one can write

(3.11)
$$\begin{aligned}\mathbb{P}(yf(x) \leq 0) &\leq \mathbb{E}_\xi \mathbb{P}\varphi_\delta(yg(x)) + \mathbb{P}(\sigma^2(c,x) \geq \gamma) \\ &\quad + \mathbb{E}\mathbb{P}_\xi(yg(x) \geq \delta, yf(x) \leq 0, \sigma^2(c;x) \leq \gamma).\end{aligned}$$

Similarly to (3.9) one can also write

(3.12)
$$\begin{aligned}\mathbb{E}_\xi \mathbb{P}_n \varphi_\delta(yg(x)) &\leq \mathbb{E}_\xi \mathbb{P}_n(yg(x) \leq 2\delta) \\ &\leq \mathbb{P}_n(yf(x) \leq 3\delta) + \mathbb{P}_n(\sigma^2(c;x) \geq \gamma) \\ &\quad + \mathbb{P}_n \mathbb{P}_\xi(yg(x) \leq 2\delta, yf(x) \geq 3\delta, \sigma^2(c;x) \leq \gamma).\end{aligned}$$

*Step* 2 (Bernstein's inequality). To bound the last terms on the right-hand sides of (3.11) and (3.12) we note that we explicitly introduced the condition on the variance of the $g_k$'s, since for a fixed $x \in \mathcal{X}$ we have $\text{Var}_\xi(g_k(x)) = \sigma^2(c; x)$. Therefore, instead of using Hoeffding's inequality as we did in the proof of Theorem 2, it is advantageous to use Bernstein's inequality, since it takes into account the information about the variance. We have

$$\mathbb{P}_\xi(yg(x) \geq \delta, yf(x) \leq 0, \sigma^2(c;x) \leq \gamma)$$

$$\leq \mathbb{P}_\xi\left(\sum_{k=1}^N (yg_k(x) - y\mathbb{E}_\xi g_k(x)) \geq N\delta \mid \text{Var}_\xi(g_1(x)) \leq \gamma\right)$$

$$\leq \exp\left(-\frac{1}{4}\min\left(\frac{N\delta^2}{\gamma}, N\delta\right)\right) = \exp\left(-\frac{1}{4}\frac{N\delta^2}{\gamma}\right),$$

since we assume that $\gamma \geq \delta$. Taking $N = 4(\gamma/\delta^2) \log n$ we get

(3.13)  $\mathbb{P}(yf(x) \leq 0) \leq \mathbb{E}_\xi \mathbb{P}\varphi_\delta(yg(x)) + \mathbb{P}(\sigma^2(c;x) \geq \gamma) + n^{-1}.$

Similarly, applying Bernstein's inequality to the last term of (3.12) yields

(3.14)  $\mathbb{E}_\xi \mathbb{P}_n \varphi_\delta(yg(x)) \leq \mathbb{P}_n(yf(x) \leq 3\delta) + \mathbb{P}_n(\sigma^2(c;x) \geq \gamma) + n^{-1}.$

*Step* 3 [Relating $\mathbb{E}_\xi \mathbb{P}\varphi_\delta(yg(x))$ to $\mathbb{E}_\xi \mathbb{P}_n \varphi_\delta(yg(x))$]. Our next goal is to relate $\mathbb{E}_\xi \mathbb{P}\varphi_\delta(yg(x))$ from the right-hand side of (3.13) to $\mathbb{E}_\xi \mathbb{P}_n \varphi_\delta(yg(x))$ from the left-hand side of (3.14).

For any realization of random variables $\xi_k^i$, the function $g(x)$ will belong to the class $\mathcal{F}_{mN}$. Convexity of the function $\phi(a,b)$ and Theorem 7 imply that for any $t > 0$ with probability at least $1 - e^{-t}$ for all $\delta \in \Delta$, $\lambda \in \mathcal{P}(\mathcal{H})$ and $f(x) = \int h(x) \, d\lambda$, and any $c \in \mathcal{C}^m(\lambda)$,

$$\phi(\mathbb{E}_\xi \mathbb{P}\varphi_\delta(yg(x)), \mathbb{E}_\xi \mathbb{P}_n \varphi_\delta(yg(x))) \leq \mathbb{E}_\xi \phi(\mathbb{P}\varphi_\delta(yg(x)), \mathbb{P}_n \varphi_\delta(yg(x)))$$

$$\leq K\left(\frac{VmN}{n}\log\frac{n}{\delta} + \frac{t}{n}\right).$$

The fact that $\phi(a,b)$ is decreasing in $b$ and increasing in $a$ combined with (3.13) and (3.14) [recall that $N = 4(\gamma/\delta^2)\log n$] implies that

$$\phi(\mathbb{P}(yf(x) \leq 0) - \mathbb{P}(\sigma^2(c;x) \geq \gamma) - n^{-1},$$

$$\mathbb{P}_n(yf(x) \leq 3\delta) + \mathbb{P}_n(\sigma^2(c;x) \geq \gamma) + n^{-1})$$

$$\leq K\left(\frac{Vm\gamma}{n\delta^2}\log^2\frac{n}{\delta} + \frac{t}{n}\right).$$

Solving the last inequality for $\mathbb{P}(yf(x) \leq 0)$ one can get that with probability at least $1 - e^{-t}$ for all $\delta \in \Delta$, any $\gamma \geq \delta$, for any $\lambda \in \mathcal{P}(\mathcal{H})$ and $f(x) =$

$\int h(x)\,d\lambda$, and any $c \in \mathcal{C}^m(\lambda)$,

$$\mathbb{P}(yf(x) \leq 0) \leq K\bigg(\mathbb{P}_n(yf(x) \leq 3\delta) + \mathbb{P}_n(\sigma^2(c;x) \geq \gamma)$$
(3.15)
$$+ \mathbb{P}(\sigma^2(c;x) \geq \gamma) + \frac{Vm\gamma}{\delta^2}\log^2\frac{n}{\delta} + \frac{t}{n}\bigg).$$

*Step* 4 [Bounding $\mathbb{P}(\sigma^2(c;x) \geq \gamma)$]. It remains to estimate $\mathbb{P}(\sigma^2(c;x) \geq \gamma)$. This is done very similarly to steps 1–3 above. Let us generate two independent sequences $\xi_k^{j,1}$ and $\xi_k^{j,2}$ as above and consider

$$\sigma_N^2(c;x) = \frac{1}{2N}\sum_{k=1}^{N}\bigg(\sum_{j=1}^{m}\alpha_j(\xi_k^{j,1} - \xi_k^{j,2})\bigg)^2 = \frac{1}{N}\sum_{k=1}^{N}\xi_k(x),$$

where

(3.16)
$$\xi_k(x) = \tfrac{1}{2}\bigg(\sum_{j=1}^{m}\alpha_j(\xi_k^{j,1} - \xi_k^{j,2})\bigg)^2.$$

Let us make a specific choice of functions $\varphi_\gamma$. For each $\gamma \in \Delta$ we set $\varphi_\gamma$ to be $\varphi_\gamma(s) = 0$ for $s \leq 2\gamma$, $\varphi_\gamma(s) = 1$ for $s \geq 3\gamma$ and linear on $[2\gamma, 3\gamma]$. One can write

$$\mathbb{P}(\sigma^2(c;x) \geq 4\gamma) = \mathbb{E}_\xi\mathbb{P}(\sigma^2(c;x) \geq 4\gamma, \sigma_N^2(c;x) \geq 3\gamma)$$
$$+ \mathbb{E}_\xi\mathbb{P}(\sigma^2(c;x) \geq 4\gamma, \sigma_N^2(c;x) \leq 3\gamma)$$
(3.17)
$$\leq \mathbb{E}_\xi\mathbb{P}\varphi_\gamma(\sigma_N^2(c;x))$$
$$+ \mathbb{E}\mathbb{P}_\xi(\sigma_N^2(c;x) \leq 3\gamma, \sigma^2(c;x) \geq 4\gamma).$$

Similarly, one can write

$$\mathbb{E}_\xi\mathbb{P}_n\varphi_\gamma(\sigma_N^2(c;x)) \leq \mathbb{E}_\xi\mathbb{P}_n(\sigma_N^2(c;x) \geq 2\gamma)$$
(3.18)
$$\leq \mathbb{P}_n(\sigma^2(c;x) \geq \gamma)$$
$$+ \mathbb{P}_n\mathbb{P}_\xi(\sigma_N^2(c;x) \geq 2\gamma, \sigma^2(c;x) \leq \gamma).$$

Next we will show that there exists a large enough absolute constant $K > 0$ such that

(3.19)
$$\mathbb{P}_\xi(\sigma_N^2(c;x) \geq 2\gamma, \sigma^2(c;x) \leq \gamma) \leq \exp\bigg(-\frac{N\gamma}{K}\bigg)$$

and

(3.20)
$$\mathbb{P}_\xi(\sigma_N^2(c;x) \leq 3\gamma, \sigma^2(c;x) \geq 4\gamma) \leq \exp\bigg(-\frac{N\gamma}{K}\bigg).$$

First of all, let us notice that $\sigma_N^2(c;x) = N^{-1} \sum_{i=1}^{N} \xi_k(x)$, where $\xi_k$ are i.i.d. random variables defined in (3.16) and $\mathbb{E}_\xi \xi_k(x) = \sigma^2(c;x)$. Moreover, since $\xi_k^{j,1}, \xi_k^{j,2} \in \mathcal{H}$, we have $|\xi_k^{j,1}(x) - \xi_k^{j,2}(x)| \leq 2$ and $|\xi_k(x)| \leq 2$. Finally, the variance

$$\operatorname{Var}_\xi(\xi_1) \leq \mathbb{E}_\xi \xi_1^2 \leq 2\mathbb{E}_\xi \xi_1 = 2\sigma^2(c;x).$$

Hence, Bernstein's inequality implies that

$$\mathbb{P}_\xi(\sigma_N^2(c;x) - \sigma^2(c;x) \leq 2\sqrt{\sigma^2(c;x)\gamma/K} + 8\gamma/(3K)) \geq 1 - \exp\left(-\frac{N\gamma}{K}\right)$$

and

$$\mathbb{P}_\xi(\sigma^2(c;x) - \sigma_N^2(c;x) \leq 2\sqrt{\sigma^2(c;x)\gamma/K} + 8\gamma/(3K)) \geq 1 - \exp\left(-\frac{N\gamma}{K}\right).$$

It is now easy to check that for large enough $K > 0$, given $\sigma^2(c;x) \leq \gamma$, the first inequality will imply $\sigma_N^2(c;x) \leq 2\gamma$ [with probability at least $1 - \exp(-\frac{N\gamma}{K})$], thus proving (3.19) and, given $\sigma_N^2(c;x) \leq 3\gamma$, the second inequality will similarly imply $\sigma^2(c;x) \leq 4\gamma$, thus proving (3.20).

If in (3.19) and (3.20) we set $N = K\gamma^{-1} \log n$, then with this choice of $N$ one can rewrite (3.17) and (3.18) as

(3.21) $$\mathbb{P}(\sigma^2(c;x) \geq 4\gamma) \leq \mathbb{E}_\xi \mathbb{P} \varphi_\gamma(\sigma_N^2(c;x)) + n^{-1}$$

and

(3.22) $$\mathbb{E}_\xi \mathbb{P}_n \varphi_\gamma(\sigma_N^2(c;x)) \leq \mathbb{P}_n(\sigma^2(c;x) \geq \gamma) + n^{-1}.$$

For any realization of $\xi_k^{j,1}, \xi_k^{j,2}$, the function $\sigma_N^2(c;x)$ belongs to the class

$$\mathcal{F}_{N,m} = \left\{ \frac{1}{N} \sum_{k=1}^{N} \left( \sum_{j=1}^{m} \alpha_j (h_k^{j,1} - h_k^{j,2}) \right)^2 : h_k^{j,1}, h_k^{j,2} \in \mathcal{H}, \alpha_j \geq 0, \sum_{j=1}^{m} \alpha_j = 1 \right\}.$$

Since the class $\mathcal{H}$ satisfies condition (2.2), it is easy to show (see, e.g., [21] for a similar computation) that the uniform covering numbers of $\mathcal{F}_{N,m}$ can be bounded by

$$\log N(\mathcal{F}_{N,m}, u) \leq KVNm \log \frac{2}{u}, \qquad 0 < u \leq 1.$$

The rest of the argument is similar to the above. Convexity of the function $\phi(a,b)$ and Theorem 7 imply that for any $t > 0$ with probability at least $1 - e^{-t}$ for all $\gamma \in \Delta$, $\lambda \in \mathcal{P}(\mathcal{H})$ and any $c \in \mathcal{C}^m(\lambda)$,

$$\phi(\mathbb{E}_\xi \mathbb{P} \varphi_\gamma(\sigma_N^2(c;x)), \mathbb{E}_\xi \mathbb{P}_n \varphi_\gamma(\sigma_N^2(c;x)))$$
$$\leq \mathbb{E}_\xi \phi(\mathbb{P} \varphi_\gamma(\sigma_N^2(c;x)), \mathbb{P}_n \varphi_\gamma(\sigma_N^2(c;x)))$$
$$\leq K\left(\frac{VmN}{n} \log \frac{n}{\delta} + \frac{t}{n}\right).$$

The fact that $\phi(a,b)$ is decreasing in $b$ and increasing in $a$ combined with (3.21) and (3.22) (recall that $N = K \log n/\gamma$) implies that

$$\phi(\mathbb{P}(\sigma^2(c;x) \geq 4\gamma) - n^{-1}, \mathbb{P}_n(\sigma^2(c;x) \geq \gamma) + n^{-1}) \leq K\left(\frac{Vm}{n\gamma}\log^2\frac{n}{\delta} + \frac{t}{n}\right).$$

Solving the last inequality for $\mathbb{P}(\sigma^2(c;x) \geq 4\gamma)$ we get that with probability at least $1 - e^{-t}$ for any $\gamma \in \Delta$, for any $\lambda \in \mathcal{P}(\mathcal{H})$ and any $c \in \mathcal{C}^m(\lambda)$,

$$\mathbb{P}(\sigma^2(c;x) \geq 4\gamma) \leq K\left(\mathbb{P}_n(\sigma^2(c;x) \geq \gamma) + \frac{Vm}{n\gamma}\log^2\frac{n}{\delta} + \frac{t}{n}\right).$$

Finally, we combine this with (3.15) and notice that since we assume that $\gamma \geq \delta$,

$$\frac{Vm}{\gamma}\log^2\frac{n}{\delta} \leq \frac{Vm\gamma}{\delta^2}\log^2\frac{n}{\delta}.$$

Thus, with probability at least $1 - e^{-t}$ for any $\delta \in \Delta$, any $\delta \leq \gamma \in \Delta$ for any $\lambda \in \mathcal{P}(\mathcal{H})$ and any $c \in \mathcal{C}^m(\lambda)$,

(3.23)
$$\begin{aligned}\mathbb{P}(yf(x) \leq 0) \\ \leq K\left(\mathbb{P}_n(yf(x) \leq 3\delta) + \mathbb{P}_n(\sigma^2(c;x) \geq \gamma/4) + \frac{Vm\gamma}{\delta^2}\log^2\frac{n}{\delta} + \frac{t}{n}\right).\end{aligned}$$

Using the union bound one can show that with a larger constant $K$ this inequality holds for all $m \geq 1$, also with probability at least $1 - e^{-t}$. Finally, to obtain the statement of Theorem 4, we need to make the change of variables $3\delta \to \delta$, $\gamma/4 \to \gamma$, and, in order to preserve the condition $\gamma \geq \delta$, we notice that from the very beginning we could have assumed that $\gamma \geq 12\delta$ and then deal with the case of $\gamma \in [\delta, 12\delta]$ by increasing the value of $K$. □

We turn now to the proof of Theorem 5. It will be based on several facts.
First of all, we need a slight modification of Theorem 2 in [8].
Let $\mathcal{F}$ be a class of functions $f$ from $\mathcal{X}$ into $[0,1]$. We define the Rademacher process $R_n(f)$, $f \in \mathcal{F}$, as

$$R_n(f) := n^{-1}\sum_{i=1}^{n}\varepsilon_i f(X_i),$$

where $\{\varepsilon_i\}$ is a Rademacher sequence $[\mathbf{Pr}(\varepsilon_i = 1) = \mathbf{Pr}(\varepsilon_i = -1) = 1/2]$ independent of $\{X_i\}$. Denote also

$$R_n(\mathcal{F}) := \sup_{f \in \mathcal{F}}|R_n(f)|.$$

THEOREM 8. *Suppose that for all $t > 0$ with probability at least $1 - e^{-t}$*

$$\mathbb{E}_\varepsilon \sup_{f \in \mathcal{F}, \mathbb{P}_n f \leq r} |R_n(f)| \leq C(\phi_n(\sqrt{r}) + \delta_n(t)) \qquad \text{where } r > 0,$$

*$\phi_n$ is a nondecreasing concave possibly data-dependent function with $\phi_n(0) = 0$, $\delta_n(t) \geq \frac{t}{n}$ and $C > 0$ is a constant. Let $\hat{r}_n$ be the largest solution of the equation $\phi_n(\sqrt{r}) = r$. Then, there exists $K > 0$ such that with probability at least $1 - e^{-t}$ for all $f \in \mathcal{F}$*

$$\mathbb{P}f \leq K\left(\mathbb{P}_n f + \hat{r}_n + \delta_n\left(\frac{t + \log \log n}{n}\right)\right).$$

Next we need the following bound on the expected sup-norm of the Rademacher process. Let

$$D_{\mathbb{P}_n,2}(\mathcal{F}) := \sup_{f,g \in \mathcal{F}} d_{\mathbb{P}_n,2}(f,g)$$

denote the $\mathcal{L}_2(\mathbb{P}_n)$-diameter of $\mathcal{F}$.

LEMMA 1. *Let $\mathcal{F}$ be a class of measurable functions from $\mathcal{X}$ into $[0,1]$ such that $0 \in \mathcal{F}$. Then there exists a constant $K > 0$ such that for all $n \geq 1$ and $t > 0$*

$$\mathbb{E}_\varepsilon R_n(\mathcal{F}) \leq \frac{K}{\sqrt{n}}\left[\sqrt{\frac{t}{n}} H_{d_{\mathbb{P}_n,1}}^{1/2}\left(\mathcal{F}; \frac{t}{n}\right) + \int_{\sqrt{t/n}}^{D_{\mathbb{P}_n,2}(\mathcal{F})} H_{d_{\mathbb{P}_n,2}}^{1/2}(\mathcal{F}, u)\, du\right] + \frac{t}{n}.$$

PROOF. For given $t > 0$ and $n \geq 1$, there exists a map $\pi = \pi_{n,t} : \mathcal{F} \mapsto \mathcal{F}$ such that

$$\operatorname{card}(\pi \mathcal{F}) = N_{d_{\mathbb{P}_n,2}}\left(\mathcal{F}, \sqrt{\frac{t}{n}}\right) \quad \text{and} \quad d_{\mathbb{P}_n,2}(f, \pi f) \leq \sqrt{\frac{t}{n}}.$$

This implies that

$$\mathbb{E}_\varepsilon R_n(\mathcal{F}) \leq \mathbb{E}_\varepsilon \sup_{f \in \mathcal{F}} |R_n(\pi f)| + \mathbb{E}_\varepsilon \sup_{f \in \mathcal{F}} |R_n(f - \pi f)|.$$

By a standard entropy bound, we have

$$\mathbb{E}_\varepsilon \sup_{f \in \mathcal{F}} |R_n(\pi f)| \leq \frac{K}{\sqrt{n}} \int_{\sqrt{t/n}}^{D_{\mathbb{P}_n,2}(\mathcal{F})} H_{d_{\mathbb{P}_n,2}}^{1/2}(\mathcal{F}, u)\, du$$

with some constant $K > 0$. Let now $\mathcal{F}'$ be a $(t/n)$-net for $\mathcal{F}$ with respect to the metric $d_{\mathbb{P}_n,1}$. Note that, since the functions from $\mathcal{F}$ take their values in $[0,1]$,

$$d_{\mathbb{P}_n,1}(f, f') \leq \frac{t}{n} \quad \Longrightarrow \quad d_{\mathbb{P}_n,2}^2(f, f') \leq \frac{t}{n}$$

$$\Longrightarrow \quad d_{\mathbb{P}_n,2}(f', \pi f) \leq d_{\mathbb{P}_n,2}(f, f') + d_{\mathbb{P}_n,2}(f, \pi f) \leq 2\sqrt{\frac{t}{n}}.$$

Therefore, we get

$$\mathbb{E}_\varepsilon \sup_{f \in \mathcal{F}} |R_n(f - \pi f)|$$
$$\leq \mathbb{E}_\varepsilon \sup\left\{|R_n(f' - g)| : f' \in \mathcal{F}', g \in \pi\mathcal{F}, d_{\mathbb{P}_n,2}(f', \pi f) \leq 2\sqrt{\frac{t}{n}}\right\}.$$

Since

$$\operatorname{card}(\mathcal{F}' \times \pi\mathcal{F}) = N_{d_{\mathbb{P}_n,1}}\left(\mathcal{F}, \frac{t}{n}\right) N_{d_{\mathbb{P}_n,2}}\left(\mathcal{F}, \sqrt{\frac{t}{n}}\right),$$

we get, using standard bounds for the expectation of a finite maximum of a Rademacher process,

$$\mathbb{E}_\varepsilon \sup_{f \in \mathcal{F}} |R_n(f - \pi f)| \leq K\sqrt{\frac{t}{n}} \frac{1}{\sqrt{n}} \left( H_{d_{\mathbb{P}_n,1}}\left(\mathcal{F}, \frac{t}{n}\right) + H_{d_{\mathbb{P}_n,2}}\left(\mathcal{F}, \sqrt{\frac{t}{n}}\right) \right)^{1/2} + \frac{t}{n}$$

with some $K > 0$, which in view of the trivial bound

$$H_{d_{\mathbb{P}_n,2}}\left(\mathcal{F}, \sqrt{\frac{t}{n}}\right) \leq H_{d_{\mathbb{P}_n,1}}\left(\mathcal{F}, \frac{t}{n}\right)$$

implies the statement of the lemma. $\square$

Let $Q \in \mathcal{P}(\mathcal{X})$. For a set $E$ of positive numbers and a function $N : E \mapsto \mathbb{R}_+$ let

$$\mathcal{F}_{Q,p,N}^C := \{f \in \mathcal{F} : \forall \varepsilon \in E \ N_{d_{Q,p}}(f, C\varepsilon) \leq N(\varepsilon)\}.$$

LEMMA 2. *For all $\varepsilon \in E$*

$$H_{d_{Q,p}}(\mathcal{F}_{Q,p,N}^C, (2+C)\varepsilon) \leq KN(\varepsilon) \log \frac{1}{\varepsilon}$$

*with some constant $K > 0$.*

PROOF. First note that $f \in \mathcal{F}_{Q,p,N}^C$ implies that

$$\forall \varepsilon \in E \ \exists \mathcal{H}' \subset \mathcal{H} : f \in \operatorname{sconv}(\mathcal{H}')$$

and

$$N_{d_{Q,p}}(\mathcal{H}', C\varepsilon) \leq N(\varepsilon).$$

Let $f = \sum \lambda_j h_j$, $h_j \in \mathcal{H}'$ and $\sum |\lambda_j| \leq 1$. Then there exists $\bar{\mathcal{H}}' \subset \mathcal{H}'$ such that $\operatorname{card}(\bar{\mathcal{H}}') \leq N(\varepsilon)$ and for all $h \in \mathcal{H}'$ there exists $g \in \bar{\mathcal{H}}'$ such that $d_{Q,p}(h, g) \leq C\varepsilon$. Hence, one can define $\{\bar{h}'_j\} \subset \bar{\mathcal{H}}'$ such that $\max_j d_{Q,p}(h_j, \bar{h}'_j) \leq C\varepsilon$. Let now $\mathcal{H}_\varepsilon$ denote a minimal $\varepsilon$-net for $\mathcal{H}$ with respect to $d_{Q,p}$. Define $\bar{h}_j \in \mathcal{H}_\varepsilon$ in such a way that for all $j$, $d_{Q,p}(\bar{h}_j, \bar{h}'_j) \leq \varepsilon$, which, of course, implies

$$\max_j d_{Q,p}(h_j, \bar{h}_j) \leq (C+1)\varepsilon.$$

Clearly, we can also assume that

$$\text{card}\{\bar{h}_j\} \leq \text{card}\{\bar{h}'_j\} \leq \text{card}(\bar{\mathcal{H}}') \leq N(\varepsilon).$$

We can conclude that

$$d_{Q,p}\left(\sum_j \lambda_j h_j, \sum_j \lambda_j \bar{h}_j\right) \leq \sum_j |\lambda_j| d_{Q,p}(h_j, \bar{h}_j)$$

$$\leq \sum_j |\lambda_j| \max_j d_{Q,p}(h_j, \bar{h}_j)$$

$$\leq (C+1)\varepsilon.$$

The above argument shows that $\forall \varepsilon \in E$

$$\mathcal{F}^C_{Q,p,N} \subset [\text{sconv}_{N(\varepsilon)}(\mathcal{H}_\varepsilon)]_{(C+1)\varepsilon},$$

where $[\cdot]_\varepsilon$ denotes the $\varepsilon$-neighborhood w.r.t. the metric $d_{Q,p}$ and

$$\text{sconv}_d(\mathcal{G}) := \left\{\sum_{j=1}^d \lambda_j h_j : \sum_{j=1}^d |\lambda_j| \leq 1 \ \forall j h_j \in \mathcal{G}\right\}.$$

Using Lemma 3 in [21], we obtain that $\forall \varepsilon \in E$

$$N_{d_{Q,p}}(\mathcal{F}^C_{Q,p,N}, (2+C)\varepsilon) \leq \left(\frac{e^2 \text{card}(\mathcal{H}_\varepsilon)(N(\varepsilon) + 4\varepsilon^{-2})}{N^2(\varepsilon)}\right)^{N(\varepsilon)},$$

which immediately implies the bound. $\square$

LEMMA 3. *Suppose that $\mathcal{H}$ satisfies* (2.2). *Then there exist constants $K > 0$, $C > 0$ such that for all $t > KV(\mathcal{H}) \log n$, with probability at least $1 - e^{-t}$ for all $f \in \mathcal{F}$ and all $\varepsilon \geq \sqrt{\frac{t}{n}}$*

$$N_{d_{\mathbb{P}_n,2}}(f, C\varepsilon) \leq N_{d_{\mathbb{P},2}}(f, \varepsilon)$$

*and*

$$N_{d_{\mathbb{P},2}}(f, C\varepsilon) \leq N_{d_{\mathbb{P}_n,2}}(f, \varepsilon).$$

PROOF.  Let

$$\tilde{\mathcal{H}} := \{(h_1 - h_2)^2 : h_1, h_2 \in \mathcal{H}\}.$$

Since $-1 \leq h \leq 1$ for $h \in \mathcal{H}$ one can write
$$((h_1 - h_2)^2 - (h'_1 - h'_2)^2)^2 \leq 32((h_1 - h'_1)^2 + (h_2 - h'_2)^2).$$
Hence the uniform covering numbers of $\tilde{\mathcal{H}}$ can be estimated as
$$\sup_{Q \in \mathcal{P}(\mathcal{X})} N_{d_{Q,2}}(\tilde{H}, \varepsilon) \leq \sup_{Q \in \mathcal{P}(\mathcal{X})} N^2_{d_{Q,2}}(\mathcal{H}, \varepsilon/8) = \mathcal{O}(\varepsilon^{-4V(\mathcal{H})})$$
using (2.3). Now, applying Theorem 7 and (3.5), we get that with probability at least $1 - 2e^{-t}$, for all $h \in \tilde{\mathcal{H}}$
$$\mathbb{P}h - \mathbb{P}_n h \leq K\left(\left(\frac{(\mathbb{P}h)V \log n}{n}\right)^{1/2} + \left(\frac{(\mathbb{P}h)t}{n}\right)^{1/2}\right)$$
and
$$\mathbb{P}_n h - \mathbb{P}h \leq K\left(\left(\frac{(\mathbb{P}_n h)V \log n}{n}\right)^{1/2} + \left(\frac{(\mathbb{P}_n h)t}{n}\right)^{1/2}\right).$$
For $t \geq KV \log n$ these inequalities imply
$$\mathbb{P}h \leq K\left(\mathbb{P}_n h + \frac{t}{n}\right) \quad \text{and} \quad \mathbb{P}_n h \leq K\left(\mathbb{P}h + \frac{t}{n}\right).$$
This yields that with probability $1 - 2e^{-t}$ for all $h_1, h_2 \in \mathcal{H}$
$$d_{\mathbb{P}_n, 2}(h_1, h_2) \leq C\left[d_{\mathbb{P},2}(h_1, h_2) \vee \sqrt{\frac{t}{n}}\right]$$
and
$$d_{\mathbb{P},2}(h_1, h_2) \leq C\left[d_{\mathbb{P}_n, 2}(h_1, h_2) \vee \sqrt{\frac{t}{n}}\right].$$
Now, by the definition of $N_{d_{\mathbb{P},2}}(f, \varepsilon)$, there exists $\mathcal{H}' \subset \mathcal{H}$ such that $f \in \text{sconv}(\mathcal{H}')$ and $N_{d_{\mathbb{P},2}}(\mathcal{H}', \varepsilon) = N_{d_{\mathbb{P},2}}(f, \varepsilon)$. Hence, with probability at least $1 - 2e^{-t}$, for any $\varepsilon \geq \sqrt{\frac{t}{n}}$, we have
$$N_{d_{\mathbb{P}_n,2}}(f, C\varepsilon) \leq N_{d_{\mathbb{P}_n,2}}(\mathcal{H}', C\varepsilon) \leq N_{d_{\mathbb{P},2}}(\mathcal{H}', \varepsilon) = N_{d_{\mathbb{P},2}}(f, \varepsilon),$$
and similarly
$$N_{d_{\mathbb{P},2}}(f, C\varepsilon) \leq N_{d_{\mathbb{P}_n,2}}(f, \varepsilon),$$
which immediately implies the bound of the lemma (after a minor rescaling and changing the constants). $\square$

Let us define a sequence
$$\varepsilon_j := 2^j \sqrt{\frac{t}{n}} \quad \text{for } j \geq 0.$$

Denote $m_n(t) := \min\{j : \varepsilon_j \geq 1\}$. Let $N$ be a nonnegative nonincreasing function on $\mathbb{R}_+$ taking constant values on the intervals $(0, \varepsilon_1)$, $[\varepsilon_j, \varepsilon_{j+1})$, $j \geq 1$. Define

$$\hat{\mathcal{F}}_{\mathbb{P}_n, N} := \{f \in \mathcal{F} : N_{d_{\mathbb{P}_n, 2}}(f, \varepsilon_j) \leq N(\varepsilon_j),\ j = 0, \ldots, m_n(t)\},$$

$$\mathcal{F}_{\mathbb{P}, N} := \{f \in \mathcal{F} : N_{d_{\mathbb{P}, 2}}(f, C\varepsilon_j) \leq N(\varepsilon_j),\ j = 0, \ldots, m_n(t)\},$$

$$\tilde{\mathcal{F}}_{\mathbb{P}_n, N} := \{f \in \mathcal{F} : N_{d_{\mathbb{P}_n, 2}}(f, C^2\varepsilon_j) \leq N(\varepsilon_j),\ j = 0, \ldots, m_n(t)\}.$$

Then it follows from Lemma 3 that:

LEMMA 4.
$$\mathbf{Pr}\{\hat{\mathcal{F}}_{\mathbb{P}_n, N} \subset \mathcal{F}_{\mathbb{P}, N} \subset \tilde{\mathcal{F}}_{\mathbb{P}_n, N}\} \geq 1 - e^{-t}.$$

Let us introduce the function

$$\psi(x) := \psi_N(x) := \int_0^x \sqrt{N(\varepsilon) \log \frac{1}{\varepsilon}}\, d\varepsilon.$$

LEMMA 5. *There exists $K > 0$ such that with probability at least $1 - e^{-t}$ for all $f \in \mathcal{F}_{\mathbb{P}, N}$*

$$\mathbb{P}\{yf(x) \leq 0\} \leq K \inf_{\delta \in (0, 1]} \left[\mathbb{P}_n\{yf(x) \leq \delta\} + \varepsilon_n^\psi(\delta) + \frac{t + \log \log n}{n\delta^2}\right].$$

PROOF. We apply Lemma 1 with $t$ replaced by $(2 + C^2)^2 t/\delta^2$ to the class

$$\mathcal{G} := \{\varphi \circ f : f \in \mathcal{F}_{\mathbb{P}, N}\} \cup \{0\},$$

where $\varphi$ is the function equal to 1 for $u \leq 0$, equal to 0 for $u > \delta$ and linear in between and $(\varphi \circ f)(x, y) := \varphi(yf(x))$. This gives the bound

$$\mathbb{E}_\varepsilon \sup_{g \in \mathcal{G}, \mathbb{P}_n g \leq r} |R_n(g)|$$

$$\leq \frac{K}{\sqrt{n}} \left[\frac{2 + C^2}{\delta} \sqrt{\frac{t}{n}} H_{d_{\mathbb{P}_n, 1}}^{1/2}\left(\mathcal{G}, \frac{(2 + C^2)^2 t}{n\delta^2}\right) + \int_{(2+C^2)/\delta\sqrt{t/n}}^{(2r)^{1/2}} H_{d_{\mathbb{P}_n, 2}}^{1/2}(\mathcal{G}, u)\, du\right]$$

$$+ \frac{(2 + C^2)^2 t}{n\delta^2}.$$

Since the Lipschitz norm of $\varphi$ is $\frac{1}{\delta}$, we have

$$d_{\mathbb{P}_n, 2}(\varphi \circ f, \varphi \circ g) \leq \frac{1}{\delta} d_{\mathbb{P}_n, 2}(f, g)$$

and
$$d_{\mathbb{P}_n,1}(\varphi \circ f, \varphi \circ g) \leq \frac{1}{\delta} d_{\mathbb{P}_n,1}(f,g).$$

Therefore, we can upper bound the expression in the brackets by
$$\frac{2+C^2}{\delta}\sqrt{\frac{t}{n}}\left[H_{d_{\mathbb{P}_n,1}}^{1/2}\left(\mathcal{F}_{\mathbb{P},N}, \frac{(2+C^2)^2 t}{n\delta}\right)+1\right]$$
$$+\frac{1}{\delta}\int_{(2+C^2)\sqrt{t/n}}^{\delta(2r)^{1/2}}\left(H_{d_{\mathbb{P}_n,2}}^{1/2}(\mathcal{F}_{\mathbb{P},N}, u)+1\right)du$$

[adding 1 to the square root of the entropy is due to the definition of the class $\mathcal{G}$ which includes the function 0; we also use here the inequality $\sqrt{\log(N+1)} \leq \sqrt{\log N}+1$]. On the event $\{\mathcal{F}_{\mathbb{P},N} \subset \tilde{\mathcal{F}}_{\mathbb{P}_n,N}\}$, which according to Lemma 4 occurs with probability at least $1-e^{-t}$, we can upper bound the $\mathcal{L}_2(\mathbb{P}_n)$- and $\mathcal{L}_1(\mathbb{P}_n)$-entropies involved in the last expression by the entropies of the class $\tilde{\mathcal{F}}_{\mathbb{P}_n,N}$, which can be bounded using Lemma 2. Namely, we have, for all $f \in \tilde{\mathcal{F}}_{\mathbb{P}_n,N}$,
$$N_{d_{\mathbb{P}_n,2}}(f, C^2\varepsilon_j) \leq N(\varepsilon_j), \qquad j=0,\ldots,m_n(t),$$
which according to Lemma 2 implies that
$$H_{d_{\mathbb{P}_n,2}}(\tilde{\mathcal{F}}_{\mathbb{P},N}, (2+C^2)\varepsilon_j) \leq KN(\varepsilon_j)\log(1/\varepsilon_j).$$

Therefore, denoting $\bar{\varepsilon}_j := (2+C^2)\varepsilon_j$ and using monotonicity of the entropy, we get
$$\int_{(2+C^2)\sqrt{t/n}}^{\delta(2r)^{1/2}} H_{d_{\mathbb{P}_n,2}}^{1/2}(\tilde{\mathcal{F}}_{\mathbb{P},N}, u)\,du$$
$$\leq \sum_{j:\bar{\varepsilon}_j \leq \delta(2r)^{1/2}}(\bar{\varepsilon}_{j+1}-\bar{\varepsilon}_j)H_{d_{\mathbb{P}_n,2}}^{1/2}(\tilde{\mathcal{F}}_{\mathbb{P},N},\bar{\varepsilon}_j)$$
$$\leq K\sum_{j:(2+C^2)\varepsilon_j \leq \delta(2r)^{1/2}}(2+C^2)(\varepsilon_{j+1}-\varepsilon_j)\sqrt{N(\varepsilon_j)\log(1/\varepsilon_j)}$$
$$\leq K\int_{\sqrt{t/(2n)}}^{2\sqrt{2}\delta\sqrt{r}}\sqrt{N(u)|\log u|}\,du.$$

Note also that since the class $\mathcal{H}$ consists of functions taking values in $\{-1,1\}$, for any probability measure $Q$ we have $d_{Q,2}^2(h_1,h_2)=2d_{Q,1}(h_1,h_2)$, which implies that $N_{d_{Q,2}}(f,\sqrt{2\varepsilon})=N_{d_{Q,1}}(f,\varepsilon)$. Thus,
$$\forall f \in \tilde{\mathcal{F}}_{\mathbb{P}_n,N} \quad N_{d_{\mathbb{P}_n,2}}(f, C^2\varepsilon_0) \leq N(\varepsilon_0)$$
$$\implies \forall f \in \tilde{\mathcal{F}}_{\mathbb{P}_n,N} \quad N_{d_{\mathbb{P}_n,1}}(f, C^4\varepsilon_0^2/2) \leq N(\varepsilon_0).$$

Since $\varepsilon_0 = \sqrt{\frac{t}{n}}$, this, in view of Lemma 2, yields the bound

$$H_{d_{\mathbb{P}_n,1}}\left(\tilde{\mathcal{F}}_{\mathbb{P}_n,N}, (2+C^4/2)\sqrt{\frac{t}{n}}\right) \leq KN\left(\sqrt{\frac{t}{n}}\right)\log\sqrt{\frac{n}{t}}.$$

Collecting the above bounds gives on the event $\{\mathcal{F}_{\mathbb{P},N} \subset \tilde{\mathcal{F}}_{\mathbb{P}_n,N}\}$

$$\left[\frac{2+C^2}{\delta}\sqrt{\frac{t}{n}}\left(H^{1/2}_{d_{\mathbb{P}_n,1}}\left(\mathcal{F}_{\mathbb{P},N}, \frac{(2+C^2)^2 t}{n\delta}\right) + 1\right)\right.$$

$$\left. + \frac{1}{\delta}\int_{(2+C^2)\sqrt{t/n}}^{\delta(2r)^{1/2}}(H^{1/2}_{d_{\mathbb{P}_n,2}}(\mathcal{F}_{P,N}, u) + 1)\,du\right]$$

$$\leq \frac{K}{\delta}\left[\sqrt{\frac{t}{2n}}\sqrt{N\left(\sqrt{\frac{t}{2n}}\right)\left|\log\sqrt{\frac{t}{2n}}\right|} + \int_{\sqrt{t/(2n)}}^{2\sqrt{2}\delta\sqrt{r}}\sqrt{N(u)|\log u|}\,du\right]$$

$$+ 2\sqrt{2}\sqrt{r},$$

which, using the fact that the function $x \mapsto \int_0^x \sqrt{N(u)|\log u|}\,du$ is concave, can be bounded by $K\phi_n(\sqrt{r})$, where

$$\phi_n(\sqrt{r}) := \phi_{n,\delta}(\sqrt{r}) := \frac{1}{\delta}\int_0^{\delta\sqrt{r}}\sqrt{N(u)|\log u|}\,du.$$

Thus, with probability at least $1 - e^{-t}$,

$$\mathbb{E}_\varepsilon \sup_{g \in \mathcal{G},\, \mathbb{P}_n g \leq r} |R_n(g)| \leq K\left(\phi_n(\sqrt{r}) + \frac{t}{n\delta^2}\right)$$

and Theorem 8 implies that also with probability at least $1 - e^{-t}$ for all $g \in \mathcal{G}$

$$\mathbb{P}g \leq K\left(\mathbb{P}_n g + \hat{r}_n + \frac{t + \log\log n}{n\delta^2}\right),$$

where $\hat{r}_n$ is the largest solution of the equation $\phi_n(\sqrt{r}) = r$, which in our case is equal to $\varepsilon_n^\psi(\delta)$. Therefore, for a fixed $\delta \in (0,1]$ with probability at least $1 - e^{-t}$ for all $f \in \mathcal{F}_{\mathbb{P},N}$

$$\mathbb{P}\{yf(x) \leq 0\} \leq \mathbb{P}(\varphi \circ f) \leq K\left(\mathbb{P}_n(\varphi \circ f) + \varepsilon_n^\psi(\delta) + \frac{t + \log\log n}{n\delta^2}\right)$$

$$\leq K\left(\mathbb{P}_n\{yf(x) \leq \delta\} + \varepsilon_n^\psi(\delta) + \frac{t + \log\log n}{n\delta^2}\right).$$

It remains to make the bound uniform in $\delta \in (0,1]$ by applying it with $\delta = \delta_j = 2^{-j}$ and $t$ replaced by $t + 2\log(j+1)$, using the union bound along with the monotonicity of the expressions involved with respect to $\delta$, and properly adjusting the value of the constant $K$. □

PROOF OF THEOREM 5. We will prove, in fact, an improved version of the result (see the remark after the statement). To simplify the notation, we remove the term $\varepsilon^{-2V/(2+V)}$ from the definition of $\hat{H}_n(f,\varepsilon)$ and the follow-up definition of $\hat{\psi}_n(f,t,\delta)$; this omission, however, does not change anything in the proof. By the condition on the class $\mathcal{H}$,

$$\sup_{Q \in \mathcal{P}(\mathcal{X})} N_{d_{Q,2}}(\mathcal{H},\varepsilon) = \mathcal{O}(\varepsilon^{-V}), \qquad \varepsilon > 0.$$

Clearly, we have

$$N_{d_{\mathbb{P}_n,2}}(f,\varepsilon) \leq \sup_{Q \in \mathcal{P}(S)} N_{d_{Q,2}}(\mathcal{H},\varepsilon), \qquad \varepsilon > 0.$$

As before, $\varepsilon_j = 2^j \sqrt{\frac{t}{n}}$ and let $J := \{j \geq 0 : \varepsilon_j < 2\}$. Denote by $\mathcal{N}$ the set of nonincreasing step functions on $\mathbb{R}_+$ with jumps only at the points $\varepsilon_j$, $j \geq 0$, and such that

$$N(\varepsilon_j) \leq K \varepsilon_j^{-V}, \qquad j \in J.$$

Assume also that, for $N \in \mathcal{N}$ and $\varepsilon \leq \varepsilon_0$, $N(\varepsilon) = N(\varepsilon_0)$. Then

$$\mathbf{Pr}\Big\{\exists f \in \mathcal{F} \, \exists \delta \in (0,1]:$$
$$\mathbb{P}\{yf(x) \leq 0\} \geq K\Big[\mathbb{P}_n\{yf(x) \leq \delta\} + \hat{\varepsilon}_n(f,t,\delta) + \frac{t + \log\log n}{n\delta^2}\Big]\Big\}$$
$$\leq \mathbb{E} \sum_{N \in \mathcal{N}} I(N_{d_{\mathbb{P}_n,2}}(f,\varepsilon_j) = N(\varepsilon_j), \, j \in J)$$
$$\times I\Big(\exists f \in \hat{\mathcal{F}}_{\mathbb{P}_n,N} \, \exists \delta \in (0,1]:$$
$$\mathbb{P}\{yf(x) \leq 0\} \geq K\Big[\mathbb{P}_n\{yf(x) \leq \delta\}$$
$$+ \varepsilon_n^{\psi_N}(\delta) + \frac{t + \log\log n}{n\delta^2}\Big]\Big) =: B,$$

where we used the facts that, on the event $\{N_{d_{\mathbb{P}_n,2}}(f,\varepsilon_j) = N(\varepsilon_j), \, j \in J\}$,

$$f \in \mathcal{F} \implies f \in \hat{\mathcal{F}}_{\mathbb{P},N}$$

and also we have on the same event $\hat{\psi}_n(f,t,u) \leq \psi_N(u)$, $u \geq 0$, which yields

$$\hat{\varepsilon}_n(f,t,\delta) \leq \varepsilon_n^{\psi_N}(\delta).$$

According to Lemma 4, for all $N \in \mathcal{N}$, $\hat{\mathcal{F}}_{\mathbb{P}_n,N} \subset \mathcal{F}_{\mathbb{P},N}$ with probability at least $1 - e^{-t}$. Also, by simple combinatorics,

$$\operatorname{card}(\mathcal{N}) \leq \prod_{j \in J} K \left(\frac{1}{\varepsilon_j}\right)^V.$$

Therefore, we can use Lemma 5 and further bound $B$ by

$$\sum_{N \in \mathcal{N}} \mathbb{E} I(N_{d_{\mathbb{P}_n},2}(f, \varepsilon_j) = N(\varepsilon_j), j \in J)$$

$$\times I\bigg(\exists f \in \hat{\mathcal{F}}_{\mathbb{P}_n,N} \; \exists \delta \in (0,1]:$$

$$\mathbb{P}\{yf(x) \leq 0\} \geq K\bigg[\mathbb{P}_n\{yf(x) \leq \delta\} + \varepsilon_n^{\psi_N}(\delta) + \frac{t + \log \log n}{n\delta^2}\bigg]\bigg)$$

$$\leq \prod_{j \in J} K\left(\frac{1}{\varepsilon_j}\right)^V \bigg[\sup_{N \in \mathcal{N}} \mathbf{Pr}\bigg\{\exists f \in \mathcal{F}_{\mathbb{P},N} \; \exists \delta \in (0,1]:$$

$$\mathbb{P}\{yf(x) \leq 0\}$$

$$\geq K\bigg[\mathbb{P}_n\{yf(x) \leq \delta\}$$

$$+ \varepsilon_n^{\psi_N}(\delta) + \frac{t + \log \log n}{n\delta^2}\bigg]\bigg\} + e^{-t}\bigg]$$

$$\leq 2 \exp\bigg\{-t + \sum_{j \in J}\bigg(V \log \frac{1}{\varepsilon_j} + \log K\bigg)\bigg\}$$

$$\leq 2 \exp\bigg\{-t + C \log^2 \frac{n}{t} + \log 2\bigg\},$$

which implies the bound of the theorem (subject to adjusting the constants). □

**4. Concluding remarks.** We have developed several new complexity measures of functions from the convex hull of a given base class and proved adaptive margin type bounds on the generalization error of ensemble classifiers in terms of these complexities. The complexities are based on measuring sparsity of the weights of a convex combination and clustering of the base functions involved in it. Hopefully, they can provide some insights to the developers of classification algorithms about the relative importance of various parameters influencing the performance of classifiers. It might be possible to combine several types of bounds discussed in the paper into a bound that takes into account different complexity characteristics, but our goal here is not to develop "the Mother of All Bounds," but rather to explore several possible approaches to the problem.

The results of the paper suggest that it might be of interest to study experimentally the statistical properties of base classifiers in ensembles output by classification algorithms (in particular, their clustering properties) in connection with generalization ability of the algorithms. (Some preliminary results in this direction for *AdaBoost* and other classification algorithms with real and simulated data can be found in [20] and more results are in [1].) Another interesting line of research might be related to proving that boosting type algorithms do output combined classifiers with a certain degree of clustering of base classifiers in the ensemble and a certain degree of sparsity of their weights. (The results of [30] show that the sparsity of the coefficients indeed takes place in the case of support vector machines.)

Our main goal has been to develop margin-type bounds on generalization error in terms of sparsity and clustering, but the complexities we introduced might be of interest in some other problems, for instance, in studying convergence rates of classification algorithms to the Bayes risk. Recent results on consistency [15, 22, 33, 34] and convergence rates [6, 7] of boosting-type algorithms suggest that some regularization of the algorithms (either by early stopping, or by penalization) might be needed in order to achieve reasonable convergence rates. However, the precise form of this regularization is still an open question and it depends crucially on which complexity measures are used to take into account the sparsity and the clustering properties of the algorithms. Some of the complexities discussed in the paper might be used as penalties, especially, the complexities based on the notion of variance of a convex combination (this is also computationally attractive). Another area where these complexities might be very useful is the problem of optimal aggregation of estimators in regression or classification (see [3, 31]).

It should be emphasized that the complexities of convex combinations we have introduced are by no means the only possible, but they are on the other hand very typical, representing some features of functions in the convex hull that are of importance in classification.

**Acknowledgment.** We would like to thank the anonymous referees for doing a great job that led to the improvement of this paper.

COMPLEXITY AND GENERALIZATION BOUND 43[25] PANCHENKO, D. (2002). Some extensions of an inequality of Vapnik and Chervonenkis. *Electron. Comm. Probab.* **7** 55–65. MR1887174

[26] PANCHENKO, D. (2003). Symmetrization approach to concentration inequalities for empirical processes. *Ann. Probab.* **31** 2068–2081. MR2016612

[27] PISIER, G. (1981). Remarques sur un résultat non publié de B. Maurey. In *Seminar on Functional Analysis, 1980–1981, École Polytechnic, Palaiseau.* Exp. no. V, 13 pp. MR659306

[28] SCHAPIRE, R., FREUND, Y., BARTLETT, P. and LEE, W. S. (1998). Boosting the margin: A new explanation for the effectiveness of voting methods. *Ann. Statist.* **26** 1651–1686. MR1673273

[29] SCHAPIRE, R. and SINGER, Y. (1999). Improved boosting algorithms using confidence-rated predictions. *Machine Learning* **37** 297–336. MR1811573

[30] STEINWART, I. (2004). Sparseness of support vector machines. *J. Mach. Learn. Res.* **4** 1071–1105. MR2125346

[31] TSYBAKOV, A. (2003). Optimal rates of aggregation. *Lecture Notes in Artificial Intelligence.* To appear.

[32] VAN DER VAART, A. W. and WELLNER, J. A. (1996). *Weak Convergence and Empirical Processes. With Applications to Statistics.* Springer, New York. MR1385671

[33] ZHANG, T. (2004). Statistical behavior and consistency of classification methods based on convex risk minimization. *Ann. Statist.* **32** 56–85. MR2051051

[34] ZHANG, T. and YU, B. (2005). Boosting with early stopping: Convergence and consistency. *Ann. Statist.* **33** 1538–1579.
| | |
|---|---|
| DEPARTMENT OF MATHEMATICS AND STATISTICS | DEPARTMENT OF MATHEMATICS |
| UNIVERSITY OF NEW MEXICO | MASSACHUSETTS INSTITUTE OF TECHNOLOGY |
| ALBUQUERQUE, NEW MEXICO 87131-1141 | CAMBRIDGE, MASSACHUSETTS 02139-4307 |
| USA | USA |
| E-MAIL: vlad@math.unm.edu | E-MAIL: panchenk@math.mit.edu |